\documentclass[final]{siamart1116}

\usepackage{enumerate}

\renewcommand{\hat}{\widehat}
\renewcommand{\epsilon}{\varepsilon}

\renewcommand{\l}{\left}
\renewcommand{\r}{\right}

\newcommand{\mbb}{\mathbb}

\newcommand{\rtaw}{\rightarrow}

\newcommand*{\defeq}{\mathrel{\vcenter{\baselineskip0.5ex \lineskiplimit0pt
                     \hbox{\scriptsize.}\hbox{\scriptsize.}}}%
                     =}

\renewcommand{\phi}{\varphi}

\newcommand{\tbf}[1]{\textbf{#1}}
\newcommand{\txt}[1]{\text{#1}}

\newcommand{\bsym}[1]{\boldsymbol{#1}}
\DeclareMathOperator*{\argmin}{arg\,min}
\DeclareMathOperator*{\argmax}{arg\,max}

\newcommand{\notetoself}[1]{\text{}}

\usepackage{lipsum}
\usepackage{amsfonts}
\usepackage{graphicx}
\usepackage{epstopdf}
\usepackage{algorithmic}
\ifpdf
  \DeclareGraphicsExtensions{.eps,.pdf,.png,.jpg}
\else
  \DeclareGraphicsExtensions{.eps}
\fi

\numberwithin{theorem}{section}

\newcommand{\TheTitle}{A Splitting Method For Overcoming the Curse of Dimensionality in Hamilton-Jacobi Equations Arising from Nonlinear Optimal Control and Differential Games with Applications to Trajectory Generation}
\newcommand{\TheShortTitle}{A Splitting Method to Compute Hamilton-Jacobi Equations}
\newcommand{\TheAuthors}{A. T. Lin, Y. T. Chow, and S. J. Osher}

\headers{\TheShortTitle}{\TheAuthors}

\title{{\TheTitle}\thanks{\funding{This work was partially supported by DOE-SC0013838 and NGA HM0210410003}}}

\author{
  Alex Tong Lin\thanks{Department of Mathematics, University of California, Los Angeles, 520 Portola Plaza, Math
Sciences Building 6363, Los Angeles, CA 90095
    (\email{atlin@math.ucla.edu}, \url{www.math.ucla.edu/\~atlin}), (\email{ytchow@math.ucla.edu}, \url{www.math.ucla.edu/\~ytchow})  (\email{sjo@math.ucla.edu}, \url{www.math.ucla.edu/\~sjo}).}
  \and
  Yat Tin Chow\footnotemark[2]
  \and
  Stanley J. Osher\footnotemark[2]
}

\usepackage{amsopn}
\DeclareMathOperator{\diag}{diag}

\ifpdf
\hypersetup{
  pdftitle={\TheTitle},
  pdfauthor={\TheAuthors}
}
\fi

\begin{document}

\maketitle

\begin{abstract}
  Recent observations have been made that bridge splitting methods arising from optimization, to the Hopf and Lax formulas for Hamilton-Jacobi Equations with Hamiltonians $H(p)$. This has produced extremely fast algorithms in computing solutions of these PDEs. More recent observations were made in generalizing the Hopf and Lax formulas to state-and-time-dependent cases $H(x,p,t)$. In this article, we apply a new splitting method based on the Primal Dual Hybrid Gradient algorithm (a.k.a. Chambolle-Pock) to nonlinear optimal control and differential games problems, based on techniques from the derivation of the new Hopf and Lax formulas, which allow us to compute solutions at points $(x,t)$ directly, i.e. without the use of grids in space. This algorithm also allows us to create trajectories directly. Thus we are able to lift the curse of dimensionality a bit, and therefore compute solutions in much higher dimensions than before. And in our numerical experiments, we actually observe that our computations scale polynomially in time. Furthermore, this new algorithm is embarrassingly parallelizable.
\end{abstract}

%

\section{Introduction}

	Hamilton-Jacobi Equations (HJE) are crucial in solving and analyzing problems arising from optimal control, differential games, dynamical systems, calculus of variations, quantum mechanics, and the list goes on \cite{EvaSou83, OshFed03}.
	
	Most methods to compute HJE use grids and finite-difference discretization. Some of these methods use ENO/WENO-type methods \cite{OshShu91}, and others use Dijkstra-type methods \cite{Dij59} such as fast marching \cite{411258} and fast sweeping \cite{TCOZ03}. But due to their use of grids, they suffer from the curse of dimensionality, i.e. they do not scale well with increases in dimension in the space variable, i.e. they generally scale exponentially.
	
	In past years, there has been an effort to mitigate the effects of dimensionality on computations of HJE. Some recent attempts to solve Hamilton-Jacobi equations use methods from low rank tensor representations \cite{7040310}, or methods based on alternating least squares \cite{7759553}, or methods by sparse grids \cite{Kang2017}, or methods using pseudospectral \cite{10.1007/978-3-540-45056-6_21} and iterative methods \cite{2017arXiv170204400K}. There have also been attempts to mitigate the curse of dimensionality which have been motivated by reachability \cite{DBLP:journals/corr/abs-1709-07523, MITCHELL2012108}. In this work, we examine and advertize the effectiveness of splitting to solve Hamilton-Jacobi equations and to directly compute optimal trajectories.
	
	We note that splitting for optimal control problems was used by \cite{6422363} (2013), where they applied it to cost functionals with a quadratic and convex term. In terms of Hamilton-Jacobi equations, Kirchner et al. \cite{KHDO18} (2018) have effectively applied PDHG \cite{ZhuChan08, doi:10.1137/09076934X} (a.k.a. Chambolle-Pock \cite{Chambolle2011}) to Hamilton-Jacobi equations arising from linear optimal control problems. They applied splitting to the Hopf formula to compute HJE for bounded input, high-dimensional linear control problems. Another main feature of their methods is they are able to directly generate optimal trajectories by making use of the the closed-form solution to linear ODEs. See also previous work by Kirchner et al. \cite{8231132} where they apply the Hopf formula to a differential games problems, which resulted in complex ``teaming" behavior even under linearized pursuit-evasion models.
	
	In this current paper, we have worked in parallel with the above authors and have also applied splitting to Hamilton-Jacobi equations arising from nonlinear problems. Our volunteered method has some nice properties: (1) relatively quick computations of solutions in high dimensions (see \cref{sec:quadcopter}, although one can easily extend to 100 dimensions for example, and also see \cref{subsubsec:dimscale} where we observe a linear relationship between computation time and dimension), especially when we include parallelization, the method is embarrassingly parallelizable \cite{Herlihy:2008:AMP:1734069}, (2) the ability to \emph{directly generate} optimal trajectories of the optimal control/differential games problems, (3) the ability to compute problems with non-linear ODE dynamics, (4) the ability to compute solutions for nonconvex Hamiltonians, (5) the ease of parallelization of our algorithm to compute solutions to HJE, i.e. each core can use the algorithm to compute the solution at a point, so given $N$ cores we can compute solutions of the HJE at $N$ points simultaneously, and (6) the ease of parallelization to directly compute trajectories, i.e. in our discretization of the time, we can parallelize by assigning each computational core a point in the time discretization.

	Our work lies in using the techniques used to derive the Generalized Hopf and Lax formulas introduced by Y.T. Chow, J. Darbon, S. Osher, and W. Yin \cite{CDOYApr17}, which generalize to the state-and-time-dependent cases (note in the literature that the classical Lax formula is sometimes called the Hopf-Lax formula). See also previous work from the same authors, \cite{Chow2016, Chow2017}, and also \cite{JDarbon15, JO16} where they provide fast algorithms under convexity assumptions. To perform the optimization, we use a new splitting method that is based on the Primal Dual Hybrid Gradient (PDHG) method (a.k.a. Chambolle-Pock), which we \emph{conjecture} to both converge to a local minimum for most well-behaved problems, and which we \emph{conjecture} to also approximate the solution. To do this, we discretize the optimal control problem and the differential games problem in time, a technique inspired by \cite{6422363} and \cite{CDOYApr17}. This new splitting method has been experimentally seen (\cref{sec:examples}) to be faster than the using coordinate-descent in most cases, which the authors in \cite{CDOYApr17} use to compute the solutions.
	
%
	
	As far as the authors know, the use of splitting as applied to minimax differential games problems, mainly on the state-and-time dependent equation \eqref{eq:discgenlaxhopfdg} and \eqref{eq:discgenlaxhopfdghopf}, is new. In this case, we seem to be able to compute HJE with nonconvex Hamiltonians and nonconvex initial data (\cref{subsubsec:isaacscc}), although they do have the structure of being convex-concave.

The paper is organized as follows:
\begin{itemize}
	\item \cref{sec:hjeocdg} Gives brief overviews of Hamilton-Jacobi Equations and its intimate connections to optimal control \cref{subsec:hjeoc} and differential games \cref{subsec:hjedg}.
	\item \cref{sec:split} Gives a brief overview of splitting methods from optimization, focusing on ADMM \cref{subsec:ADMM} and PDHG \cref{subsec:PDHG}.
	\item \cref{sec:genlaxhopf} Presents the generalized Lax and Hopf formulas for optimal control and differential games that were conjectured by \cite{CDOYApr17}. We also go through its discretization in \cref{sec:discgenlaxhopf}, which is the basis of our algorithm.
	\item \cref{sec:algo} \textbf{Presents the main algorithms.}
	\item \cref{sec:examples} Presents various computational examples.
	\item \cref{sec:disc} Ends with a brief conclusion, and a discussion on future work.
	\item \cref{appa:sec:practut} Gives a more in-depth explanation on how to use the algorithms.
\end{itemize}

\newpage

\section{Hamilton-Jacobi Equations and Its Connection to Optimal Control and Differential Games}
\label{sec:hjeocdg}

\subsection{Hamilton-Jacobi Equations and Optimal Control}
\label{subsec:hjeoc}

	Most of our exposition on optimal control will follow \cite{EvansOptBook}, and also \cite{evans10} (Chapter 10).

	The goal of optimal control theory is to find a control policy that will drive a system while optimizing a criterion. 
	Given an initial point $x\in\mbb{R}^n$ and an initial time $t\in [0, T]$, where $T$ is some fixed end-point time, the system will obey an ODE
		\begin{align*}
		\l\{\begin{array}{l}
		\dot{\tbf{x}}(s) = \tbf{f}(\tbf{x}(s), \bsym{\alpha}(s), s), \quad t<s<T \\
		\tbf{x}(t) = x
		\end{array}\r.
		\end{align*}
	where $\tbf{f}:(\mbb{R}^n \times A \times \mbb{R}) \rtaw \mbb{R}$, where $A\subseteq \mbb{R}^m$. We call $\tbf{x}$ the \emph{state}, and $\bsym{\alpha}$ the \emph{control}. And the functional we want to optimize is $J_{x,t}:\mathcal{A} \rtaw \mbb{R}$ where
		\begin{align}\label{eq:Jval}
		J_{x,t}[\bsym{\alpha}] \defeq g(\tbf{x}(T)) + \int_t^T L(\tbf{x}(s), \bsym{\alpha}(s), s)\,ds.
		\end{align}
	and where $\mathcal{A} \defeq \{\bsym{\alpha}:[t,T]\rtaw A\}$ is some \emph{admissable control} set, and $g:\mbb{R}^n\rtaw \mbb{R}$ and $L:(\mbb{R}^n \times A \times \mbb{R}) \rtaw \mbb{R}$. We can either \emph{minimize} the above functional, in which we call it a \emph{cost}, or we can \emph{maximize} it, in which we call it a \emph{payoff}. For our exposition, we will choose to minimize $J_{x,t}[\cdot]$, so it will be a cost. Then we define the \emph{value function}
		\begin{align}\label{eq:valfunc}
		\phi(x,t) = \min_{\alpha(\cdot)\in\mathcal{A}} J_{x,t}[\bsym{\alpha}].
		\end{align}
	Under some mild conditions on $\tbf{f}$, $g$, and $L$, this value function will satisfy the \emph{terminal-valued} Hamilton-Jacobi PDE (HJ PDE)
		\begin{align*}
		\l\{\begin{array}{rcl}
		\partial_t \phi(x,t) + H(x, \nabla_x \phi(x,t), t) &=& 0, \qquad (x,t) \in\mbb{R}^n \times (0,T)\\
		\phi(x,T) &=& g(x).
		\end{array}\r.
		\end{align*}
	where $H(x,p,t) = \min_{a\in A}\l\{ \l<\tbf{f}(x, a, t), p\r> + L(x,a,t)\r\}$. 
	
	To get an initial-valued PDE, can make a change of variables $t \rtaw T - t$. Or equivalently we can reformulate the optimal control problem ``backwards in time" so that we have
		\begin{align*}
		\l\{\begin{array}{l}
		\dot{\tbf{x}}(s) = \tbf{f}(\tbf{x}(s), \bsym{\alpha}(s), s), \quad 0 < s < t \\
		\tbf{x}(t) = x
		\end{array}\r.
		\end{align*}
	and 
		\begin{align*}
		J_{x,t}[\bsym{\alpha}] \defeq g(\tbf{x}(0)) + \int_0^t L(\tbf{x}(s), \bsym{\alpha}(s), s)\,ds.
		\end{align*}
	Then our $\phi(x,t) = \min_{\alpha(\cdot)\in\mathcal{A}} J_{x,t}[\bsym{\alpha}]$ will satisfy an \emph{initial-valued} HJ PDE with $H(x,p,t) = \max_{a\in A}\l\{ \l<\tbf{f}(x, a, t), p\r> - L(x,a,t)\r\}$. Note that if $\tbf{f} = a$, then this form of the Hamiltonian expresses $H$ as the convex conjugate \cite{OPT-003} of $L$.
	
	If we thinks from a physical perspective in which time moves forward, the first formulation feels more comfortable. If we comes from the fields of PDE or mathematical optimization, the latter formulation will feel more comfortable.
	
	So how does having a HJ PDE help us synthesize an optimal control? Using the first formulation as it feels more physically intuitive, we can heuristically argue that given an inital time $t\in (0,T]$ and a state $x\in\mbb{R}^n$, we consider the optimal ODE
		\begin{align*}
		\l\{\begin{array}{l}
		\dot{\tbf{x}}^{*}(s) = \tbf{f}(\tbf{x}^*(s), \bsym{\alpha}^*(s), s), \quad t<s<T \\
		\tbf{x}^*(t) = x
		\end{array}\r.
		\end{align*}
	where at each time $s\in (0,T)$, we choose the value of $\alpha^*(s)$ to be such that
		\begin{align*}
		&\l<\tbf{f}(\tbf{x}^*(s), \bsym{\alpha}^*(s), s), \partial_x \phi(\tbf{x}^*(s),s)\r> + L(\tbf{x}^*(s), \bsym{\alpha}^*(s), s) \\
		&\qquad = \min_{a\in A}\l\{ \l<\tbf{f}(\tbf{x}^*(s), a, s), \partial_x \phi(\tbf{x}^*(s)\r> + L(\tbf{x}^*(s),a,s)\r\} \\		
		&\qquad = H(\tbf{x}^*(s), \partial_x \phi(\tbf{x}^*(s),s), s)
		\end{align*}
	We call $\bsym{\alpha}^*(\cdot)$ defined in this way as the \emph{feedback control}, and this can be obtained from $\nabla_x \phi$ (see Section 10.3.3 of \cite{evans10}). This is also related to Pontryagin' Maximum Principle (Chapter 4 of \cite{EvansOptBook}).
	
	Note that in the case $H(x,p,t) = H(p)$, then we have available the (classical) Hopf and Lax formulas which are expressions for the solutions $\phi(x,t)$ of the HJ PDE:
	
	When the Hamiltonian $H(p)$ is convex and the initial date $g$ is (uniformly) Lipschitz continuous, then we have the Lax formula:
		\begin{align*}
		\phi(x,t) = \min_{y\in\mbb{R}^n} \l\{g(y) + t H^*\l(\frac{x-y}{t}\r) \r\}
		\end{align*}
	where $H^*(x) = \max_{v\in A} \l\{\l<v,x\r> - H(v) \r\}$ is the convex conjugate of $L$. 
	
	And if the initial data $g$ is (uniformly) Lipschitz continuous and convex, and $H$ is continuous, then we have the Hopf formula,
		\begin{align}\label{eq:hopf}
		\phi(x,t) = \sup_{y\in\mbb{R}^n} \l\{-g^*(y) + \l<y,x\r> - tH(y)\r\} = \l( g^*(y) + tH(y) \r)^*
		\end{align}
	where $g^*$ is the convex conjugate of $g$. Note the last equality implies the solution is convex in $x$. We note again that the argument minimum of the above expression is in-fact $\nabla_x \phi(x,t)$.

\subsection{Hamilton-Jacobi Equations and Differential Games}
\label{subsec:hjedg}	
	
	Our exposition of differential games will follow \cite{EvansOptBook} (Chapter 6), but also see \cite{IsaacsDiffGames, YongDiffGames}. In the field of differential games, we restrict our exposition to \emph{two-person, zero-sum differential games}. Let an initial point $x\in\mbb{R}^n$ and an initial time $t\in[0,T]$ be given, where $T$ is some fixed endpoint time. A two-person, zero-sum differential game will have the dynamics,
		\begin{align*}
		\l\{\begin{array}{l}
		\dot{\tbf{x}}(s) = \tbf{f}(\tbf{x}(s), \bsym{\alpha}(s), \bsym{\beta}(s), s), \quad t<s<T \\
		\tbf{x}(t) = x
		\end{array}\r.
		\end{align*}
	where $\tbf{f}:(\mbb{R}^n \times A \times B \times \mbb{R}) \rtaw \mbb{R}$, and where $A\subseteq \mbb{R}^m$ and $B\subseteq \mbb{R}^\ell$. The control $\bsym{\alpha}$ is the control for player $I$, and the control $\bsym{\beta}$ is the control for player $II$. The functional will be $J_{x,t}:{A}(t)\times {B}(t) \rtaw \mbb{R}$ where
		\begin{align}\label{eq:dgobjfunc}
		J_{x,t}[\bsym{\alpha}, \bsym{\beta}] \defeq g(\tbf{x}(T)) + \int_t^T L(\tbf{x}(s), \bsym{\alpha}(s), \bsym{\beta}(s), s)\,ds.
		\end{align}
	and where ${A}(t) \defeq \{\bsym{\alpha}:[t,T]\rtaw A\}$ and ${B}(t) \defeq \{\bsym{\beta}:[t,T]\rtaw B\}$ are \emph{admissable control} sets, and $g:\mbb{R}^n\rtaw \mbb{R}$ and $L:(\mbb{R}^n \times A \times B \times \mbb{R}) \rtaw \mbb{R}$.
	
	In order to model that at each time, neither player has knowledge of the other's future moves, we use a concept of strategy that was used by Varaiya \cite{varaiya1972} as well as Elliot and Kalton \cite{elliott1972}. This idea allows us to model that each player will select a control in response to all possible controls the opponent can select. 
	
	A \emph{strategy} for player $I$ is a mapping $\Phi:B(t) \rtaw A(t)$ such that for all times $s\in [t,T]$
		\begin{align*}
		\tau \in [t,s], \quad \bsym{\beta}(\tau) \equiv \hat{\bsym{\beta}}(\tau) \qquad \text{implies} \qquad \Phi[\bsym{\beta}](\tau) \equiv \Phi[\hat{\bsym{\beta}}](\tau)
		\end{align*}
	The $\Phi[\bsym{\beta}]$ models player $I$'s response to player $II$ selecting $\bsym{\beta}$. We similarly define a strategy $\Psi:A(t) \rtaw B(t)$ for player $II$:
		\begin{align*}
		\tau \in [t,s], \quad \bsym{\alpha}(\tau) \equiv \hat{\bsym{\alpha}}(\tau) \qquad \text{implies} \qquad \Psi[\bsym{\alpha}](\tau) \equiv \Psi[\hat{\bsym{\alpha}}](\tau)
		\end{align*}
	and $\Psi[\bsym{\alpha}]$ models player $II$'s response to player $I$ selecting $\bsym{\alpha}$.
	
	Letting $\mathcal{A}(t)$ and $\mathcal{B}(t)$ be the set of strategies for player $I$ and player $II$, respectively, then we define the \emph{lower value function} as
		\begin{align}\label{eq:lowervalfunc}
		\phi^-(x,t) = \inf_{\Psi[\cdot]\in\mathcal{B}(t)} \sup_{\bsym{\alpha}(\cdot)\in {A}(t)} J_{x,t}[\bsym{\alpha}, \Psi[\bsym{\alpha}]]
		\end{align}
	and the \emph{upper value function} as
		\begin{align}\label{eq:uppervalfunc}
		\phi^+(x,t) = \sup_{\Phi[\cdot]\in\mathcal{A}(t)} \inf_{\bsym{\beta}(\cdot)\in B(t)} J_{x,t}[\Phi[\bsym{\beta}], \bsym{\beta}].
		\end{align}
	Note that we always have $\phi^-(x,t) \le \phi^+(x,t)$ for all $x\in\mbb{R}^n$ and $t\in [0,T]$. For a proof, see \cite{EvaSou83}.
	
	These value functions satisfy the \emph{terminal-valued} HJ PDEs
		\begin{align*}
		\l\{\begin{array}{l}
		\partial_t \phi^-(x,t) + \max_{a\in A}\min_{b\in B} \l\{ \l< \tbf{f}(x,a,b,t), \nabla_x \phi^-(x,t)\r> + L(x,a,b,t) \r\} = 0 \\
		\phi^-(x,T) = g(x)
		\end{array}\r.
		\end{align*}
	and
		\begin{align*}
		\l\{\begin{array}{l}
		\partial_t \phi^+(x,t) + \min_{b\in B}\max_{a\in A} \l\{ \l< \tbf{f}(x,a,b,t), \nabla_x \phi^+(x,t)\r> + L(x,a,b,t) \r\} = 0 \\
		\phi^+(x,T) = g(x)
		\end{array}\r.
		\end{align*}
	where we have the \emph{lower PDE Hamiltonian}
		\begin{align*}
		H^-(x,p,t) = \max_{a\in A}\min_{b\in B} \l\{ \l< \tbf{f}(x,a,b,t), p\r> + L(x,a,b,t) \r\}
		\end{align*}
	and the \emph{upper PDE Hamiltonian}
		\begin{align*}
		H^+(x,p,t) = \min_{b\in B}\max_{a\in A} \l\{ \l< \tbf{f}(x,a,b,t), p\r> + L(x,a,b,t) \r\}
		\end{align*}
	In general, we have
		\begin{align*}
		\max_{a\in A}\min_{b\in B} \l\{ \l< \tbf{f}(x,a,b,t), p\r> + L(x,a,b,t) \r\} \le \min_{b\in B}\max_{a\in A} \l\{ \l< \tbf{f}(x,a,b,t), p\r> + L(x,a,b,t) \r\}
		\end{align*}
	and in most cases the inequality is strict, and thus the lower and upper value functions are different. But when the above is an equality, then the game is said to satisfy the \emph{minimax conditions}, also called \emph{Isaac's condition}, and we have $\phi^- = \phi^+$, and we say the game has \emph{value}.
	
	Our examples will focus on differential games which satisfy the minimax condition, and will thus have value.
	
	In differential games, we usually run into nonconvex Hamiltonians, so it is the Hopf formula (\ref{eq:hopf}) that is used the most. 

\newpage
\section{Splitting Algorithms from Optimization}
\label{sec:split}
	Here we review a couple of splitting algorithms from optimization.
	
\subsection{ADMM (Alternating Method of Multipliers)}
\label{subsec:ADMM}
	
	ADMM \cite{DBLP:journals/ftml/BoydPCPE11}, which is also known as Split-Bregman \cite{doi:10.1137/080725891}, is an optimization method to solve problems of the following form:
		\begin{align*}
		\min_{x,z\in X}\quad  f(x) + g(z)\\
		\text{subject to} \quad Ax + Bz = c
		\end{align*}
	where $X$ is a finite-dimensional real vector space equipped with an inner product $\l<\cdot, \cdot\r>$, and $f:X \rtaw \mbb{R}$ and $g:X \rtaw \mbb{R}$ are proper, convex, lower semicontinuous functions. We also have that $A$ and $B$ are continuous linear operators (e.g. matrices), with $c$ a fixed element in $X$. Now we form the \emph{augmented Lagrangian} of the above problem:
		\begin{align*}
		L_\rho(x,z,y) = f(x) + g(z) + \l<y, Ax + Bz - c \r> + \frac{\rho}{2} \|Ax + Bz - c\|_2^2
		\end{align*}
	Then we alternately minimize:
		\begin{align*}
		\l\{\begin{array}{ll}
		x^{k+1} &= \argmin_x L_\rho(x, z^k, y^k) \\ \\
		z^{k+1} &= \argmin_z L_\rho(x^{k+1}, z, y^k) \\ \\
		y^{k+1} &= y^k + \rho(Ax^{k+1} + Bz^{k+1} - c)
		\end{array}\r.
		\end{align*}
	where in the last step we update the dual variable. Note that the $\argmin$ expressions are frequently precisely the \emph{proximal operator} \cite{OPT-003} of a (not necessarily convex) function. The proximal operator is defined as: Given $f:\mbb{R}^n \rtaw \mbb{R}$ a proper l.s.c. function, not necessarily convex, then,
		\begin{align}\label{eq:prox}
		(I + \lambda \partial f)^{-1}(v) \defeq \argmin_x\l\{ f(x) + \frac{1}{2\lambda} \|x-v\|_2^2 \r\}
		\end{align}
	The proximal of $f$ with step-size $\lambda$ is also denoted $\text{prox}_{\lambda f}(\cdot)$.

\subsection{PDHG (Primal-Dual Hybrid Gradient)}
\label{subsec:PDHG}
	
	The PDHG algorithm \cite{ZhuChan08, doi:10.1137/09076934X}, which also goes by the name Chambolle-Pock \cite{Chambolle2011}, attempts to solve problems of the form
		\begin{align*}
		\min_{x\in X}\quad  f(Ax) + g(x)
		\end{align*}	
	where we make similar assumptions on $X$, $f$, $g$, and $A$ as we did for ADMM. PDHG takes the Lagrangian dual formulation of the above problem and seeks to find a saddle point of the following problem:
		\begin{align*}
		\min_{x\in X} \max_{y\in Y} \l<Ax, y\r> + g(x) - f^*(y)
		\end{align*}
	where $f^*(y) = \sup_{x\in X} \l\{ \l<x,y \r> - f(y)\r\}$ is the \emph{convex conjugate} of $f$. PDHG is also an alternating minimization technique that makes use of proximal operators. The updates are:
		\begin{align*}
		\l\{\begin{array}{ll}
		y^{k+1} &= (I + \sigma \partial f^*)^{-1}(y^k + \sigma A \bar{x}^k) \\ \\
		x^{k+1} &= (I + \tau \partial g)^{-1}(x^k - \sigma A^* y^{k+1}) \\ \\
		\bar{x}^{k+1} &= x^{k+1} + \theta (x^{k+1} - x^k).
		\end{array}\r.
		\end{align*}
	where $\sigma, \tau > 0$ are such that $\sigma \tau \|A\|^2 < 1$, and $\theta \in [0,1]$, although $\theta = 1$ seems to work best in practice.
	
\newpage
\section{The Generalized Lax and Hopf formulas}
\label{sec:genlaxhopf}
	
	A recent result by Y.T. Chow, J. Darbon. S. Osher, and W. Yin \cite{CDOYApr17} gives a conjectured generalization to the Lax and Hopf formulas. Given a Hamilton-Jacobi Equation,
		\begin{align*}
		\l\{\begin{array}{rcl}
		\partial_t \phi + H(x, \nabla_x \phi(x,t), t) &=& 0, \quad \text{in } \mbb{R}^d \times (0,\infty), \\
		\phi(x,0) &=& g(x).
		\end{array}\r.
		\end{align*}
	we have that when $H(x,p,t)$ is smooth, and convex with respect to $p$, and possibly under some more mild conditions, we have
		\begin{align*}
		\phi(x,t) &= \min_{v\in\mbb{R}^d} \l\{ g(\textbf{x}(0)) + \int_0^t  \tbf{p}(s) \cdot \nabla_p H(\tbf{x}(s), \tbf{p}(s), s) - H(\tbf{x}(s),\tbf{p}(s), s) \,ds \r\}\\
		\qquad &\text{where } \l\{\begin{array}{ll}
		\dot{\tbf{x}}(s) &= \nabla_p H(\tbf{x}(s), \tbf{p}(s)) \\
		\dot{\tbf{p}}(s) &= - \nabla_x H(\tbf{x}(s), \tbf{p}(s)) \\
		\tbf{x}(t) &= x \\
		\tbf{p}(t) &= v
		\end{array}\r.
		\end{align*}
	where $\tbf{x}$ and $\tbf{p}$ are the characteristics of the PDE. The expression in the bottom braces are ODEs which $\textbf{x}(\cdot)$ and $\textbf{p}(\cdot)$ satisfy.
	
	And when we move the convexity onto $g$, i.e. when $H(x,p,t)$ is smooth and $g$ is convex, then
		\begin{align*}
		\hspace{-1cm}\varphi(x,t) &= \sup_{v\in\mbb{R}^d} \l\{ - g^\star(\tbf{p}(0)) + x\cdot v + \int_0^t  \tbf{x}(s) \cdot \nabla_x H(\tbf{x}(s), \tbf{p}(s), s) - H(\tbf{x}(s),\tbf{p}(s), s) \,ds \r\}\\
		\qquad &\text{where } \l\{\begin{array}{ll}
		\dot{\tbf{x}}(s) &= \nabla_p H(\tbf{x}(s), \tbf{p}(s)) \\
		\dot{\tbf{p}}(s) &= - \nabla_x H(\tbf{x}(s), \tbf{p}(s)) \\
		\tbf{x}(t) &= x \\
		\tbf{p}(t) &= v 
		\end{array}\r.
		\end{align*}	
		
	Chow, Darbon, Osher, and Yin used coordinate descent with multiple initial guesses to perform the optimization. They do this by first making an initial guess for $v\in\mbb{R}^d$, then they compute the ODEs, and then compute the value of the objective, i.e. the first lines of the two formulas. Then they re-adjust one coordinate $v\in\mbb{R}^d$ and repeat. Details can be found in their paper \cite{CDOYApr17}.

\subsection{Discretizing the Generalized Lax and Hopf Formulas for Optimal Control}
\label{sec:discgenlaxhopf}

	In order to derive the generalized Lax and Hopf formulas, we can first discretize the value function of the optimal control problem \eqref{eq:Jval} and \eqref{eq:valfunc}. Before we begin we note we are merely making formal calculations, much in the spirit of E. Hopf in his seminal paper where he derived the classical Hopf formula \cite{10.2307/24901315}. This is the procedure followed in \cite{CDOYApr17}: We have the value function equals
		\begin{align*}
		\phi(x,t) = \min_{\textbf{x}(\cdot), \tbf{u}(\cdot)} \l\{ g(\tbf{x}(0)) + \int_0^t L(\textbf{x}(s), \tbf{u}(s), s)\,ds \r\}
		\end{align*}
	where $\tbf{x}(\cdot)$ and $\tbf{u}(s)$ satisfy the ODE
		\begin{align*}
		\l\{\begin{array}{l}
		\dot{\tbf{x}}(s) = \tbf{f}(\tbf{x}(s), \tbf{u}(s), s), \quad 0 < s < t \\
		\tbf{x}(t) = x
		\end{array}\r.
		\end{align*}
	Two notes: (1) we are formulating our optimal control problem ``backwards in time" so that we end up with an \emph{initial-valued} HJ PDE, and (2) here $(x,t)$ are fixed points where we want to compute the HJE. 
	
	We discretize the time domain such that
		\begin{align*}
		0 < s_1 < s_2 < \cdots < s_N = t,
		\end{align*}
	and we set $x_j = \tbf{x}(s_j)$ and $u_j = \tbf{u}(s_j)$. Note in our numerical examples, we make a uniform discretization of the time domain.
	
	Now we use the backward Euler discretization of the ODE and set $x_N = x$ to obtain the optimization problem
		\begin{align*}
		\min_{\{x_j\}_{j=0}^{N}, \{u_j\}_{j=1}^{N}} \l\{ g(x_0) + \delta \sum_{j=1}^{N} L(x_j, u_j, s_j) \mid \l\{x_{j} - x_{j-1} = \delta f(x_{j}, u_{j}, s_{j})\r\}_{j=1}^N, x_N = x \r\}
		\end{align*}
	As usual in constrained optimization problems, we compute the Lagrangian function (i.e. Lagrange multipliers) to get:
		\begin{align}\label{eq:lmultformin}
		g(x_0) + \delta \sum_{j=1}^{N} L(x_j, u_j, s_j) + \sum_{j=1}^{N} \l<p_{j}, x_{j} - x_{j-1} - \delta f(x_{j}, u_{j}, s_{j})\r> + \l<p_N, x - x_N\r>
		\end{align}
	Note that the constraint $x_N = x$ is trivially unneeded in the Lagrangian function. Then we minimize over $\{x_j\}_{j=0}^{N}$ and $\{u_j\}_{j=1}^{N}$, while maximizing over $\{p_j\}_{j=1}^{N}$ to obtain the expression,
		\begin{align*}
		\max_{\{p_j\}_{j=1}^{N}} \min_{\{x_j\}_{j=0}^{N-1}} \min_{\{u_j\}_{j=1}^{N}} \l\{ g(x_0) + \delta \sum_{j=1}^{N} L(x_j, u_j, s_j) + \sum_{j=1}^{N} \l<p_{j}, x_{j} - x_{j-1} - \delta f(x_j, u_j, s_j)\r> + \l<p_N, x - x_N\r> \r\}.
		\end{align*}
	After moving the minimum over $\{u_j\}_{j=1}^{N}$ inside, we get
		\begin{align*}
		&\max_{\{p_j\}_{j=1}^{N}} \min_{\{x_k\}_{j=0}^{N}} \l\{ g(x_0) + \sum_{j=1}^{N} \l<p_{j}, x_{j} - x_{j-1}\r> + \l<p_N, x - x_N\r> + \delta \sum_{j=1}^{N} \min_{u_j} \l\{L(x_j, u_j, s_j) - \l<p_j, f(x_j, u_j, s_j)\r>\r\} \r\} \\
		&= \max_{\{p_j\}_{j=1}^{N}} \min_{\{x_k\}_{j=0}^{N}} \l\{ g(x_0) + \sum_{j=1}^{N} \l<p_{j}, x_{j} - x_{j-1}\r> + \l<p_N, x - x_N\r> - \delta \sum_{j=1}^{N} H(x_j, p_j, s_j) \r\}
		\end{align*}
	After the above step, we have now been able to remove a numerical optimization in $u$ by using the definition of the Hamiltonian. This considerably simplifies the problem, and reduces the dimensionality of the optimization.		
		
	We note that we need $p_N = 0$ in order for the maximization/minimization to not be infinite. And thus, we can remove the minimization with respect to $x^N$ and we get
		\begin{align}\label{eq:disclaxBE}
		\phi(x,t) \approx \max_{\{p_j\}_{j=1}^{N}} \min_{\{x_j\}_{j=0}^{N}} \l\{ g(x_0) + \sum_{j=1}^{N} \l<p_{j}, x_{j} - x_{j-1}\r>  - \delta \sum_{j=1}^{N} H(x_j, p_j, s_j) \r\}
		\end{align}
	We can do a similar analysis using the forward Euler discretization to obtain,
		\begin{align}\label{eq:disclaxFE}
		\phi(x,t) \approx \max_{\{p_j\}_{j=0}^{N-1}} \min_{\{x_j\}_{j=0}^{N-1}} \l\{ g(x_0) + \sum_{j=0}^{N-1} \l<p_{j}, x_{j+1} - x_j\r>  - \delta \sum_{j=0}^{N-1} H(x_j, p_j, s_j) \r\}
		\end{align}
	but this latter expression has the disadvantage of coupling the $g$ and $H$ with respect to $x_0$, as they both depend on $x_0$. Although this could actually be an advantage as one may have $H$ acting as a regularizer to $g$.
	
	In order to obtain the discretized version of the generalized Hopf formula, we start with the Lax formula with backward Euler \eqref{eq:disclaxBE}, and use the linear term $\l<p_1,x_0\r>$ and compute the convex conjugate; the calculation goes as follows:		
		\begin{align*}
		&\max_{\{p_j\}_{j=1}^{N}} \min_{\{x_j\}_{j=0}^{N}} \l\{ g(x_0) + \sum_{j=1}^{N-1} \l<p_{j}, x_{j} - x_{j-1}\r>  - \delta \sum_{j=1}^{N} H(x_j, p_j, s_j) \r\} \\
		&=\max_{\{p_j\}_{j=1}^{N}} \min_{\{x_j\}_{j=1}^{N}} \l\{ \min_{x_0}\l\{g(x_0) - \l<p_1,x_0\r>\r\} + \l<p_1,x_1\r> + \sum_{j=2}^{N-1} \l<p_{j}, x_{j} - x_{j-1}\r>  - \delta \sum_{j=1}^{N} H(x_j, p_j, s_j) \r\} \\
		&=\max_{\{p_j\}_{j=1}^{N}} \min_{\{x_j\}_{j=1}^{N}} \l\{ -g^*(p_1) + \l<p_1,x_1\r> + \sum_{j=2}^{N-1} \l<p_{j}, x_{j} - x_{j-1}\r>  - \delta \sum_{j=1}^{N} H(x_j, p_j, s_j) \r\} \\
		&= \max_{\{p_j\}_{j=1}^N} \min_{\{x_j\}_{j=1}^N } \l\{ -g^*(p_1) + \l<p_N, x\r> + \sum_{j=1}^{N-1} \l<p_j - p_{j+1}, x_j\r>  - \delta \sum_{j=1}^{N} H(x_{j}, p_{j}, s_{j})\r\} 
		\end{align*}
	where in the last equality, we performed a summation-by-parts and also used $x_N = x$. So we have the discretized version of the Hopf formula:
		\begin{align}\label{eq:dischopf}
		\phi(x,t) \approx \max_{\{p_j\}_{j=1}^N} \min_{\{x_j\}_{j=1}^N } \l\{ -g^*(p_1) + \l<p_N, x\r> + \sum_{j=1}^{N-1} \l<p_j - p_{j+1}, x_j\r>  - \delta \sum_{j=1}^{N} H(x_{j}, p_{j}, s_{j})\r\} 
		\end{align}
	Note that it is a bit harder to perform the optimization when we approximate the ODE dynamics with forward Euler because we then must compute the convex conjugate of the sum $g(\cdot) + H(\cdot, p_0, s_0)$, which can be more complicated.

\subsection{Discretizing the Generalized Lax and Hopf Formulas for Differential Games}
\label{sec:discgenlaxhopfdg}

	Again following the procedure of \cite{CDOYApr17}, we have a conjectured generalization to the Lax and Hopf formulas for differential games, which we will discretize. Before we give the calculation we qualify that, in the spirit of E. Hopf when he computed the Hopf formula in his seminal paper \cite{10.2307/24901315}, these calculations are merely formal: 
	
	Given a two-person, zero-sum differential game with value (i.e., it satisfies the Isaacs conditions so that the minmax Hamiltonian and maxmin Hamilton are equal, see section \ref{subsec:hjedg}), with given $x \in \mbb{R}^{d_1}$ and $y\in\mbb{R}^{d_2}$, and $t\in (0,\infty)$, and with dynamics
		\begin{align*}
		\l\{ \begin{array}{l}
		\begin{pmatrix}\dot{\tbf{x}}(s) \\ \dot{\tbf{y}}(s)) \end{pmatrix}= \begin{pmatrix}\tbf{f}_1(\tbf{x}(s), \tbf{y}(s), \bsym{\alpha}(s), \bsym{\beta}(s), s) \\ \tbf{f}_2(\tbf{x}(s), \tbf{y}(s), \bsym{\alpha}(s), \bsym{\beta}(s), s) \end{pmatrix} \quad 0<s<t \\
		\begin{pmatrix} \tbf{x}(t) \\ \tbf{y}(t) \end{pmatrix} = \begin{pmatrix} x \\ y \end{pmatrix}
		\end{array}\r.
		\end{align*}
	we have that the value function satisfies
		\begin{align}\label{eq:valuefuncdg}
		\phi(x,y,t) = \inf_{\bsym{\alpha}(\cdot), \tbf{x}(\cdot)} \sup_{\bsym{\beta}(\cdot), \tbf{y}(\cdot)} \l\{ g(\tbf{x}(0), \tbf{y}(0)) + \int_0^t L(\tbf{x}(s), \tbf{y}(s), \bsym{\alpha}(s), \bsym{\beta}(s), s)\,ds\r\}
		\end{align}
	
	Now, we discretize in time and approximate the ODE with backward Euler, and a formal computation gives us,
		\begin{align*}
		\approx &\min_{\{\alpha_k\}_{k=1}^N \{x_k\}_{k=0}^N} \max_{\{\beta_k\}_{k=1}^N, \{y_k\}_{k=0}^N} \l\{ g(x_0, y_0) + \delta \sum_{k=1}^N L(x_k, y_k, \alpha_k, \beta_k, s_k) \r\}  \\
		&\text{such that} \l\{ \begin{array}{l}
		\begin{pmatrix}x_k - x_{k-1} \\ y_k - y_{k-1}\end{pmatrix} = \begin{pmatrix}f_1(x_k, y_k, \alpha_k, \beta_k, s_k) \\ f_2(x_k, y_k, \alpha_k, \beta_k, s_k)\end{pmatrix}, \quad k=1, \ldots, N \\
		\begin{pmatrix} x_N \\ y_N \end{pmatrix} = \begin{pmatrix} x \\ y \end{pmatrix}
		\end{array}\r.
		\end{align*}
	It is at this point we want to form the Lagrangian. The only trouble is that the concept of a ``Lagrangian" for minimax problems (a.k.a. saddle-point problems) has not been well-examined. But in a paper by \cite{Dvurechensky2015} (and also in \cite{Qi1995}), they have a version of a Lagrangian for minimax problems, which we apply to our problem to get, 
		\begin{align*}
		g(x_0, y_0) + \delta \sum_{j=1}^N L(x_j, y_j, \alpha_j, \beta_j, s_j) + \sum_{j=1}^N \l<p_j, x_j - x_{j-1} - \delta f_1(x_j, y_j, \alpha_j, \beta_j, s_j)\r> \\
		+ \sum_{j=1}^N \l<-q_j, y_j - y_{j-1} - \delta f_2(x_j, y_j, \alpha_j, \beta_j, s_j)\r> \\
		\end{align*}
	Then we take the $\min\max$ to obtain
		\begin{align*}
		&\min_{\{\alpha_j\}_{j=1}^N, \{x_j\}_{j=0}^N} \max_{\{\beta_j\}_{j=1}^N, \{y_j\}_{j=0}^N} \l\{ g(x_0, y_0) + \delta \sum_{j=1}^N L(x_j, y_j, \alpha_j, \beta_j, s_j) \r. \\
		& \l. \hspace{0.5in} + \sum_{j=1}^N \l<p_j, x_j - x_{j-1} - \delta f_1(x_j, y_j, \alpha_j, \beta_j, s_j)\r> + \sum_{j=1}^N \l<-q_j, y_j - y_{j-1} - \delta f_2(x_j, y_j, \alpha_j, \beta_j, s_j)\r> \r\}\\
		&=\min_{\{x_j\}_{j=0}^N} \max_{\{y_j\}_{j=0}^N} \l\{ g(x_0, y_0) + \delta \sum_{j=1}^N \min_{\alpha_j}\max_{\beta_j} \l\{L(x_j, y_j, \alpha_j, \beta_j, s_j) - \l<\begin{pmatrix} p_j \\ -q_j \end{pmatrix} , \begin{pmatrix} f_1(x_j, y_j, \alpha_j, \beta_j, s_j) \\ f_2(x_j, y_j, \alpha_j, \beta_j, s_j)\end{pmatrix}\r> \r\} \r. \\
		&\hspace{2in} \l. + \sum_{j=1}^N \l<p_j, x_j - x_{j-1}\r> + \sum_{j=1}^N \l<-q_j, y_j - y_{j-1}\r> \r\}\\
		&=\min_{\{x_j\}_{j=0}^N} \max_{\{y_j\}_{j=0}^N} \l\{ g(x_0, y_0) - \delta \sum_{j=1}^N H(x_j, y_j, p_j, -q_j, s_j) + \sum_{j=1}^N \l<p_j, x_j - x_{j-1}\r> + \sum_{j=1}^N \l<-q_j, y_j - y_{j-1}\r> \r\}\\
		\end{align*}
	Now we maximize over $\{p_j\}_{j=1}^N$ and minimize over $\{q_j\}_{j=1}^N$ to obtain,
		\begin{align*}
		\phi(x,y,t) &\approx \min_{\{q_j\}_{j=1}^N} \max_{\{p_j\}_{j=1}^N} \min_{\{x_j\}_{j=0}^N} \max_{\{y_j\}_{j=0}^N} \l\{ g(x_0, y_0) - \delta \sum_{j=1}^N H(x_j, y_j, p_j, -q_j, s_j) + \sum_{j=1}^N \l<p_j, x_j - x_{j-1}\r> \r. \\
		&  \l. \hspace{4in} + \sum_{j=1}^N \l<-q_j, y_j - y_{j-1}\r> \r\}
		\end{align*}
	and after organizing a bit, we get,
		\begin{align}\label{eq:discgenlaxhopfdg}
		\phi(x,y,t) \approx \min_{\{q_j\}_{j=1}^N} \max_{\{p_j\}_{j=1}^N} \min_{\{x_j\}_{j=0}^N} \max_{\{y_j\}_{j=0}^N} \l\{ g(x_0, y_0) + \sum_{j=1}^N \l<\begin{pmatrix}p_j \\ -q_j \end{pmatrix}, \begin{pmatrix}x_j - x_{j-1} \\ y_j - y_{j-1} \end{pmatrix}\r> - \delta \sum_{j=1}^N H(x_j, y_j, p_j, -q_j, s_j)  \r\}
		\end{align}

	Note: If we can split $g(x,y) = e(x) + h(y)$, and if $e$ is convex and $h$ is concave, then we may take advantage of $e^*$, the convex conjugate of $e$, and $h_*$, the concave conjugate of $h$ (the formula for the concave conjugate is the same as the convex-conjugate, but you change the $\sup$ to an $\inf$) \cite{doi:10.1137/1.9781611970524}, in order to have an analogous Hopf formula:
	
		
		\begin{equation}\label{eq:discgenlaxhopfdghopf}
		\begin{aligned}
		\phi(x,y,t) &\approx \min_{\{q_j\}_{j=1}^N} \max_{\{p_j\}_{j=1}^N} \min_{\{x_j\}_{j=1}^N} \max_{\{y_j\}_{j=1}^N} \l\{ -e^*(p_1) - h_*(-q_1) + \l<\begin{pmatrix} p_N \\ -q_N \end{pmatrix}, \begin{pmatrix} x \\ y \end{pmatrix} \r> + \sum_{j=1}^{N-1} \l<\begin{pmatrix}p_j - p_{j+1} \\ -(q_j - q_{j+1}) \end{pmatrix}, \begin{pmatrix}x_j \\ y_j \end{pmatrix}\r>\r. \\
		& \hspace{4.5in}\l. - \delta \sum_{j=1}^N H(x_j, y_j, p_j, -q_j, s_j)  \r\}
		\end{aligned}
		\end{equation}
		
	In some ways the function $e^*(p) + h_*(q)$ may perhaps be called the \emph{convex-concave conjugate} for the convex-concave function $g(x,y)$.
	
	\textbf{Remark:} The authors in \cite{Dvurechensky2015} state that this ``minimax Lagrangian," even in the simplest formulation given in their work, is new.
	
		
\subsection{The advantage of the Hamiltonian for optimization}
\label{subsec:hamcomm}		
	
	There is tremendous advantage in having a Hamiltonian. This is because if we want to instead perform optimization of the value function directly, we will be solving for the controls and this requires a constrained optimization technique.	
	
	The miraculous advantage of having a Hamiltonian for optimization purposes is it \emph{encodes} information from both the running cost function $L$, as well as the dynamics $\dot{\tbf{x}}(s) = \tbf{f}(\tbf{x}(s), \tbf{u}(s), s)$. And thus we are now free to perform \emph{unconstrained optimization}. But if we solve for the value function \eqref{eq:Jval} and \eqref{eq:valfunc} directly, the we need to perform constrained optimization.
	
	Another key additional advantage is that we lower the dimension of the numerical optimization by analytically minimizing over $u$, and conjuring the Hamiltonian.
	

\newpage
\section{The Main Algorithm: Splitting for Hamilton-Jacobi Equations}
\label{sec:algo}

\subsection{Splitting for HJE arising from Optimal Control}
	
	Before discussing the algorithms, we note that we do not yet have a proof of convergence nor approximation. This is currently a work-in-progress. But as shown in our numerical results \cref{sec:examples}, these algorithms seem to agree with classical methods used to solve Hamilton-Jacobi equations.
	
	Taking the Lax formula with backward Euler \eqref{eq:disclaxBE} as an expository example, we can organize our problem to look similar to a primal-dual formulation which is attacked by splitting using PDHG. We stack variables and let 
		\begin{itemize}
		\item $\tilde{x} = (x_0, x_1, \ldots, x_N)$, and similarly for $\tilde{p}$ and $\tilde{s}$,
		\item $\tilde{G}(\tilde{x}) = g(x_0)$,
		\item $\tilde{H}_\delta(\tilde{x}, \tilde{p}, \tilde{s}) = \delta \sum_{k=1}^N H(x_k, u_k, s_k)$,
		\item $D$ be the difference matrix such that $\l<\tilde{p}, D\tilde{x}\r> = \sum_{j=1}^{N-1} \l<p_j, x_j - x_{j-1}\r>$
		\end{itemize}			
	then our problem looks like:
		\begin{align*}
		\max_{\{p_j\}_{j=1}^N} \min_{\{x_j\}_{j=0}^N} \tilde{G}(\tilde{x}) + \l<\tilde{p}, D\tilde{x}\r> + \tilde{H}_\delta(\tilde{x}, \tilde{p}, \tilde{s}).
		\end{align*}
	This looks similar to the problem that is attacked by PDHG, except for a couple of differences:
		\begin{itemize}
		\item PDHG solves a saddle point problem  where the $\tilde{H}_\delta(\tilde{x}, \tilde{p}, \tilde{s})$ term does not depend on $\tilde{x}$ (nor $\tilde{s}$).
	\item In PDHG, the $\tilde{H}_\delta$ term is the convex conjugate of some function we want to minimize. But in our case, we have
		\begin{align*}
		H(x,p,s) = \max_{u}\l\{ \l<f(x,u,s),p\r> - L(x,u,s) \r\}
		\end{align*}
	and $f(x,u,s)$ does not even have to be linear. So in some ways, $H$ is a ``generalized convex conjugate."
		\end{itemize}
	But we perform an alternating minimization technique similar to PDHG and we arrive at our main algorithm for optimal control: 
	
	For the Lax with backward Euler: Given an $x\in\mbb{R}^d$ and some time $t\in (0, \infty)$, we set $\delta>0$ small and let the time-grid size be $N = t/\delta + 1$ (we set $N = t/\delta$ for Hopf). Then we randomly initialize $\tilde{x} = (x_0, x_1, \ldots, x_N)$ but let $x_N = x$, and we randomly initialize $\tilde{p} = (p_0, p_1, \ldots, p_N)$ but let $p_0 \equiv 0$ as it is not minimized over, but used for computational accounting. And we let $\tilde{z} = \tilde{x}$. Then our algorithm follows the pattern of alternating optimization with quadratic penalty:
		\begin{align*}
		\l\{\begin{array}{l}
		\tilde{p}^{k+1} = \argmax_{\tilde{p}} \l\{ \tilde{G}(\tilde{x}^k) - \tilde{H}_{\delta}(\tilde{x}^k, \tilde{p}, \tilde{s}) - \frac{1}{2\sigma}\|\tilde{p} - (\tilde{p}^k + D\tilde{z}^k) \|_2^2\r\} \\
		\tilde{x}^{k+1} = \argmin_{\tilde{x}} \l\{ \tilde{G}(\tilde{x}) - \tilde{H}_{\delta}(\tilde{x}, \tilde{p}^{k+1}, \tilde{s}) + \frac{1}{2\tau}\|\tilde{x} - (\tilde{x}^{k} - D^T\tilde{p}^{k+1})\|_2^2
		\r\} \\
		\tilde{z}^{k+1} = \tilde{x}^{k+1} + \theta (\tilde{x}^{k+1} - \tilde{x}^k).
		\end{array}\r.
		\end{align*}
	where $\sigma, \tau > 0$ are step-sizes with $\sigma \tau \|D\|^2 < 1$ and $\theta \in [0,1]$ (as suggested in \cite{Chambolle2011}). In our numerical experiments, $\theta = 1$ was frequently the best choice and also	in practice, we would change the arg max into an arg min. So we have \cref{alg:algohjeoclax}.

	\begin{algorithm}
	\caption{Splitting for HJE for Optimal Control, Lax with backward Euler}
	\label{alg:algohjeoclax}
	\begin{algorithmic}
	\STATE{Given: $x_{\text{target}}\in\mbb{R}^d$ and time $t\in (0, \infty)$.}
	\STATE{Initialize: $\delta>0$ and set $N = t/\delta+1$. And randomly initialize $\tilde{x}^0$ and $\tilde{p}^0$, but with $x_N^0 \equiv x_{\text{target}}$. And set $\tilde{z}^0 = \tilde{x}^0$. Also choose $\sigma, \tau$ such that $\sigma \tau \|D\|^2 < 1$ and $\theta \in [0,1]$.}
	\WHILE{tolerance criteria large}
	\STATE{$\tilde{p}^{k+1} = \argmin_{\tilde{p}} \l\{ -\tilde{G}(\tilde{x}^k) + \tilde{H}_{\delta}(\tilde{x}^k, \tilde{p}, \tilde{s}) + \frac{1}{2\sigma}\|\tilde{p} - (\tilde{p}^k + \sigma D\tilde{z}^k) \|_2^2\r\}$}
	\STATE{$\tilde{x}^{k+1} = \argmin_{\tilde{x}} \l\{ \tilde{G}(\tilde{x}) - \tilde{H}_{\delta}(\tilde{x}, \tilde{p}^{k+1}, \tilde{s}) + \frac{1}{2\tau}\|\tilde{x} - (\tilde{x}^{k} - \tau D^T\tilde{p}^{k+1})\|_2^2
		\r\}$}
	\STATE{$\tilde{z}^{k+1} = \tilde{x}^{k+1} + \theta (\tilde{x}^{k+1} - \tilde{x}^k)$}
	\ENDWHILE
	\STATE{$\text{fval} = g(x_0) + \sum_{j=1}^{N} \l<p_{j}, x_{j} - x_{j-1}\r>  - \delta \sum_{j=1}^{N} H(x_j, u_j, s_j)$}
	\RETURN fval
	\end{algorithmic}
	\end{algorithm}

	And a similar algorithm will be obtained when we use a forward Euler discretization \cref{eq:disclaxFE} for the ODE dynamics. We can obtain better accuracy if we average the backward Euler and forward Euler approximations for the ODEs, which is reminiscent of the trapezoidal approximation having better accuracy as it is the average of the forward and backward Euler.

	We also have the Hopf formulation: Let,
		\begin{itemize}
		\item $\tilde{G}^*(\tilde{p}) = (g^*(p_1), 0, \ldots, 0)$, and
		\item let $D$ be the difference matrix such that $\l<D\tilde{p}, \tilde{x}\r> = \l<p_N, x\r> + \sum_{k=1}^{N-1} \l<p_k - p_{k+1}, x_k\r>$ (so this one differs from \cref{alg:algohjeoclax} as it acts on $\tilde{p}$)
		\end{itemize}			
	then we have \cref{alg:algohjeochopf}.
	
	\subsubsection{When to use the Hopf formula}
	
	We make the observation that the Lax formula is suitable (i.e. converges) when we have a convex Hamiltonian in $p$ which is also bounded below in $p$ (or satisfies a coercivity condition, see \cite{evans10}; if we want a convex Hamiltonian that is not bounded in $p$, then we must use Hopf in this case.

	\begin{algorithm}
	\caption{Splitting for HJE for Optimal Control, Hopf (with backward Euler)}
	\label{alg:algohjeochopf}
	\begin{algorithmic}
	\STATE{Given: $x_{\text{target}}\in\mbb{R}^d$ and time $t\in (0, \infty)$.}
	\STATE{Initialize: $\delta>0$ and set $N = t/\delta$. And randomly initialize $\tilde{x}^0$ and $\tilde{p}^0$, but with $x_N^0 \equiv x_{\text{target}}$. And set $\tilde{z}^0 = \tilde{x}^0$. Also choose $\sigma, \tau$ such that $\sigma \tau \|D\|^2 < 1$ and $\theta \in [0,1]$.}
	\WHILE{tolerance criteria large}
	\STATE{$\tilde{p}^{k+1} = \argmin_{\tilde{p}} \l\{ \tilde{G}^*(\tilde{p}^k) + \tilde{H}_{\delta}(\tilde{x}^k, \tilde{p}, \tilde{s}) + \frac{1}{2\sigma}\|\tilde{p} - (\tilde{p}^k + \sigma D^T\tilde{z}^k) \|_2^2\r\}$}
	\STATE{$\tilde{x}^{k+1} = \argmin_{\tilde{x}} \l\{ -\tilde{G}^*(\tilde{p}) - \tilde{H}_{\delta}(\tilde{x}, \tilde{p}^{k+1}, \tilde{s}) + \frac{1}{2\tau}\|\tilde{x} - (\tilde{x}^{k} - \tau D\tilde{p}^{k+1})\|_2^2
		\r\}$}
	\STATE{$\tilde{z}^{k+1} = \tilde{x}^{k+1} + \theta (\tilde{x}^{k+1} - \tilde{x}^k)$}
	\ENDWHILE
	\STATE{$\text{fval} = -g^*(p_1) + \l<p_N, x_{\text{target}}\r> + \sum_{j=1}^{N-1} \l<p_j - p_{j+1}, x_j\r>  - \delta \sum_{j=1}^{N} H(x_{j}, p_{j}, s_{j})$}
	\RETURN fval
	\end{algorithmic}
	\end{algorithm}

%

	See \cref{subsec:remarksargminmax} on how to perform the argmin/argmax in each iteration.

\subsection{Splitting for HJE arising from Differential Games}
\label{subsec:splitforhjedg}

	
	For differential games, we use a similar algorithm to the optimal control case. We take the discretized version of the cost function \eqref{eq:valuefuncdg} and perform an alternating minimization technique inspired by PDHG, but applied to minimax problems. Using the same notation as in \cref{alg:algohjeoclax} and \cref{alg:algohjeochopf}, and the same $D$ matrix as in \cref{alg:algohjeoclax}, we have the algorithm for differential games in \cref{alg:algohjedglax}. 
	
		\begin{algorithm}
		\caption{Splitting for HJE for Differential Games, Lax}
		\label{alg:algohjedglax}
		\begin{algorithmic}
		\STATE{Given: $(x_{\text{target}},y_{\text{target}})\in\mbb{R}^d$ and time $t\in (0, \infty)$.}
		\STATE{Initialize: $\delta>0$ and set $N = t/\delta+1$. And randomly set $\tilde{x}^0$, $\tilde{y}^0$, $\tilde{p}^0$, and $\tilde{q}^0$, but with $x^0_N \equiv x_{\text{target}}$ and $y^0_N \equiv y_{\text{target}}$. And set $(\bar{\tilde{x}}^0, \bar{\tilde{y}}^0) = (\tilde{x}^0, \tilde{y}^0)$. Also choose $\sigma, \tau$ such that $\sigma \tau \|D\|^2 < 1$ and $\theta \in [0,1]$.}
		\WHILE{tolerance criteria large}
		\STATE{$\tilde{p}^{k+1} = \argmax_{\tilde{p}} \l\{ \tilde{G}(\tilde{x}^k, \tilde{y}^k) - \tilde{H}_\delta(\tilde{x}^k,\tilde{y}^k,\tilde{p},-\tilde{q}^k, \tilde{s}^k) - \frac{1}{2\sigma} \|\tilde{p} - (\tilde{p}^k + \sigma D\bar{\tilde{x}}^k)\|_2^2\r\}$}
		\STATE{$\tilde{q}^{k+1} = \argmin_{\tilde{q}} \l\{ \tilde{G}(\tilde{x}^k, \tilde{y}^k) - \tilde{H}_\delta(\tilde{x}^k,\tilde{y}^k,\tilde{p}^{k+1},-\tilde{q}, \tilde{s}^k) + \frac{1}{2\sigma} \|\tilde{q} - (\tilde{q}^k + \sigma D \bar{\tilde{y}}^k)\|_2^2\r\}$}
		\STATE{$\tilde{x}^{k+1} = \argmin_{\tilde{x}} \l\{ \tilde{G}(\tilde{x}, \tilde{y}^k) - \tilde{H}_\delta(\tilde{x},\tilde{y}^k,\tilde{p}^{k+1},-\tilde{q}^{k+1}, \tilde{s}^k) + \frac{1}{2\tau} \|\tilde{x} - (\tilde{x}^k - \tau D^T\tilde{p}^{k+1})\|_2^2\r\}$}
		\STATE{$\tilde{y}^{k+1} = \argmax_{\tilde{y}} \l\{ \tilde{G}(\tilde{x}^{k+1}, \tilde{y}) - \tilde{H}_\delta(\tilde{x}^{k+1},\tilde{y},\tilde{p}^{k+1},-\tilde{q}^{k+1}, \tilde{s}^k) - \frac{1}{2\tau} \|\tilde{y} - (\tilde{y}^k - \tau D^T \tilde{q}^{k+1})\|_2^2\r\}$}
		\STATE{$\begin{pmatrix}\bar{\tilde{x}}^{k+1} \\ \bar{\tilde{y}}^{k+1}\end{pmatrix} = \begin{pmatrix}\tilde{x}^{k+1} \\ \tilde{y}^{k+1}\end{pmatrix} + \theta \l( \begin{pmatrix}\tilde{x}^{k+1} \\ \tilde{y}^{k+1}\end{pmatrix} - \begin{pmatrix}\tilde{x}^{k} \\ \tilde{y}^{k}\end{pmatrix} \r).$}
		\ENDWHILE
		\STATE{$\text{fval} = g(x_0, y_0) + \sum_{j=1}^N \l<\begin{pmatrix}p_j \\ -q_j \end{pmatrix}, \begin{pmatrix}x_j - x_{j-1} \\ y_j - y_{j-1} \end{pmatrix}\r> - \delta \sum_{j=1}^N H(x_j, y_j, p_j, -q_j, s_j)$}
		\RETURN fval
		\end{algorithmic}
		\end{algorithm}
		
	If we have $\tilde{G}(x,y) = \tilde{E}(x) + \tilde{H}(y)$ where $E$ is convex and $H$ is concave, then we may make use of convex-conjugates and concave-conjugates (see \eqref{eq:discgenlaxhopfdghopf}) to obtain an analogous Hopf formula as in \cref{alg:algohjeochopf}, but for differential games. Here $D$ is the same difference matrix as in the Hopf case, \cref{alg:algohjeochopf}. This is \cref{alg:algohjedghopf}.
		
		\begin{algorithm}
		\caption{Splitting for HJE for Differential Games, Hopf (for separable convex-concave initial conditions)}
		\label{alg:algohjedghopf}
		\begin{algorithmic}
		\STATE{Given: $(x_{\text{target}},y_{\text{target}})\in\mbb{R}^d$ and time $t\in (0, \infty)$.}
		\STATE{Initialize: $\delta>0$ and set $N = t/\delta$. And randomly set $\tilde{x}^0$, $\tilde{y}^0$, $\tilde{p}^0$, and $\tilde{q}^0$, but with $x^0_N \equiv x_{\text{target}}$ and $y^0_N \equiv y_{\text{target}}$. And set $(\bar{\tilde{x}}^0, \bar{\tilde{y}}^0) = (\tilde{x}^0, \tilde{y}^0)$. Also choose $\sigma, \tau$ such that $\sigma \tau \|D\|^2 < 1$ and $\theta \in [0,1]$.}
		\WHILE{tolerance criteria large}
		\STATE{$\tilde{p}^{k+1} = \argmax_{\tilde{p}} \l\{ -\tilde{E}^*(p_1) - \tilde{H}_*(-q_1) - \tilde{H}_\delta(\tilde{x}^k,\tilde{y}^k,\tilde{p},-\tilde{q}^k, \tilde{s}^k) - \frac{1}{2\sigma} \|\tilde{p} - (\tilde{p}^k + \sigma D^T\bar{\tilde{x}}^k)\|_2^2\r\}$}
		\STATE{$\tilde{q}^{k+1} = \argmin_{\tilde{q}} \l\{ -\tilde{E}^*(p_1) - \tilde{H}_*(-q_1) - \tilde{H}_\delta(\tilde{x}^k,\tilde{y}^k,\tilde{p}^{k+1},-\tilde{q}, \tilde{s}^k) + \frac{1}{2\sigma} \|\tilde{q} - (\tilde{q}^k + \sigma D^T \bar{\tilde{y}}^k)\|_2^2\r\}$}
		\STATE{$\tilde{x}^{k+1} = \argmin_{\tilde{x}} \l\{ -\tilde{E}^*(p_1) - \tilde{H}_*(-q_1) - \tilde{H}_\delta(\tilde{x},\tilde{y}^k,\tilde{p}^{k+1},-\tilde{q}^{k+1}, \tilde{s}^k) + \frac{1}{2\tau} \|\tilde{x} - (\tilde{x}^k - \tau D\tilde{p}^{k+1})\|_2^2\r\}$}
		\STATE{$\tilde{y}^{k+1} = \argmax_{\tilde{y}} \l\{ -\tilde{E}^*(p_1) - \tilde{H}_*(-q_1) - \tilde{H}_\delta(\tilde{x}^{k+1},\tilde{y},\tilde{p}^{k+1},-\tilde{q}^{k+1}, \tilde{s}^k) - \frac{1}{2\tau} \|\tilde{y} - (\tilde{y}^k - \tau D\tilde{q}^{k+1})\|_2^2\r\}$}
		\STATE{$\begin{pmatrix}\bar{\tilde{x}}^{k+1} \\ \bar{\tilde{y}}^{k+1}\end{pmatrix} = \begin{pmatrix}\tilde{x}^{k+1} \\ \tilde{y}^{k+1}\end{pmatrix} + \theta \l( \begin{pmatrix}\tilde{x}^{k+1} \\ \tilde{y}^{k+1}\end{pmatrix} - \begin{pmatrix}\tilde{x}^{k} \\ \tilde{y}^{k}\end{pmatrix} \r).$}
		\ENDWHILE
		\STATE{$\text{fval} = -e^*(p_1) - h_*(-q_1) + \l<\begin{pmatrix} p_N \\ -q_N \end{pmatrix}, \begin{pmatrix} x \\ y \end{pmatrix} \r> + \sum_{j=1}^{N-1} \l<\begin{pmatrix}p_j - p_{j+1} \\ -(q_j - q_{j+1}) \end{pmatrix}, \begin{pmatrix}x_j \\ y_j \end{pmatrix}\r> - \delta \sum_{j=1}^N H(x_j, y_j, p_j, -q_j, s_j)$}
		\RETURN fval
		\end{algorithmic}
		\end{algorithm}
		
\subsection{Remarks on how to perform the argmin/argmax in each iteration}
\label{subsec:remarksargminmax}
	In each iteration for the above algorithms, we have an optimization problem when updating $\tilde{x}^{k+1}$ ($\tilde{y}^{k+1}$) or $\tilde{p}^{k+1}$ ($\tilde{q}^{k+1}$). In some of our experiments, the optimization turned into a closed-form proximal expression (mainly when updating $\tilde{p}^{k+1}/\tilde{q}^{k+1}$, or we were able to make use of one step of gradient descent of the objective (mainly when updating $\tilde{x}^{k+1}$ or $\tilde{y}^{k+1}$ or when $\tilde{G} = g$ is involved).
	As an example, one way to update the Lax formula for optimal control (\cref{alg:algohjeoclax}) using a backward Euler discretization is
		\begin{align*}
		\l\{\begin{array}{ll}
		\tilde{p}^{k+1} &= \text{prox}_{\sigma \tilde{H}_\delta(\bar{\tilde{x}}^k, \cdot)} (\tilde{p}^k + \sigma D \bar{\tilde{x}}^k) \\ \\
		{\tilde{x}}^{k+1} &= {\tilde{x}}^k - \tau D^T {\tilde{p}}^{k+1} - \tau \nabla_{\tilde{x}} \tilde{G}(\tilde{x}^k) + \tau \nabla_{{\tilde{x}}} \tilde{H}_\delta({\tilde{x}}^k, {\tilde{p}}^{k+1}, {\tilde{s}}) \\ \\
		\bar{\tilde{x}}^{k+1} &= \tilde{x}^{k+1} + \theta (\tilde{x}^{k+1} - \tilde{x}^k).
		\end{array}\r.
		\end{align*}
	and one way to update the Hopf formula for optimal control (\cref{alg:algohjeochopf}) is,
		\begin{align*}
		\l\{\begin{array}{ll}
		\tilde{p}^{k+1} &= \text{prox}_{\sigma \tilde{H}_\delta(\bar{\tilde{x}}^k, \cdot)} (\tilde{p}^k + \sigma D \bar{\tilde{x}}^k - \sigma \nabla g^*(p_1)) \\ \\
		\tilde{x}^{k+1} &= \tilde{x}^k - \tau D^T \tilde{p}^{k+1} + \tau \nabla_{\tilde{x}}\tilde{H}_\delta(\tilde{x}^k, \tilde{p}^{k+1}) \\ \\
		\bar{\tilde{x}}^{k+1} &= \tilde{x}^{k+1} + \theta (\tilde{x}^{k+1} - \tilde{x}^k).
		\end{array}\r.
		\end{align*}
	where we see the update for $\tilde{p}^{k+1}$ is a \emph{proximal gradient} update \cite{OPT-003} (Section 4.3).
	
	There are many combinations we can use, but the intuition is to take a gradient step on the smooth part, and a proximal step on the non-smooth part. And if we have a sum of a smooth part and a non-smooth part, we can mix a gradient step with a proximal step (i.e. a proximal gradient step) as we have done above.
	
\subsection{Computation of characteristic curves/optimal trajectories}

	A benefit from the new algorithm is that alongside computing solutions at each point, it also allows us to directly compute trajectories/characteristic curves of the Hamilton-Jacobi equation at each point. In fact, our algorithm is a hybrid of the direct collocation method (a direct method) \cite{vonStryk1993} and Pontryagin's Maximum Principle (an indirect method) \cite{doi:10.1137/1037043}. 
	
	 We give some examples of characteristic curves/optimal trajectories in our numerical results (\cref{sec:examples}).

\subsection{The advantage of splitting over coordinate descent}
	The advantage of these methods over coordinate descent is not only its speed, but it also does not seem to require as many multiple initial guesses for nonconvex optimization. And in our numerical experiments in \cref{sec:examples}, we only used a single initial guess in all our examples. And for most examples in our experiments in \cref{sec:examples}, it only requires one guess. 
	
	It also the advantage that one can apply the method to non-smooth problems, as opposed to coordinate-descent, where one takes numerical gradients.
	
	And practically, splitting is more straightforward to implement than coordinate descent, where we would require divided differences to numerically compute the gradients, and we also have available to us the multitude of splitting techniques from the optimization literature, such as ADMM and Douglas-Rachford splitting to name a few.
	
\subsection{A remark on the connection between Hamilton-Jacobi equations and optimization, and the implications on future optimization methods}
\label{subsec:remarkhjeopt}

	The relationships between Hamilton-Jacobi equations and optimization have been noted in the literature \cite{doi:10.1137/S0363012998345366, doi:10.1137/S0363012998345378}. More recently, there has been a connection between deep learning optimization and HJE \cite{DBLP:journals/corr/ChaudhariOOSC17}. More concretely, there is a straight-forward connection between Hamilton-Jacobi equations and the proximal operator (which can be interpreted as implicit gradient descent):
	
	Given a function $f(\cdot)$, the proximal operator of $f$ is
		\begin{align*}
		(I + \lambda \partial f)^{-1}(v) \defeq \argmin_x\l\{ f(x) + \frac{1}{2\lambda} \|x-v\|_2^2 \r\}
		\end{align*}
	If we change the $\argmin_x$ into a $\min_x$, then we get the familiar Lax (a.k.a. Hopf-Lax) formula for the Hamilton-Jacobi equation with $H(x) = \frac{1}{2} \|x\|_2^2$, and for $t=\lambda$. So in this way, the Generalized Lax and Hopf formulas can be viewed as a generalization of the proximal operator.
	
	Also, the proximal operator, i.e. the $\argmin_x$ operator, is featured heavily in our algorithms. In classical PDHG, the primal variable $x$ and the dual variable $p$ are decoupled. But for our algorithms, the coupling is in the form of a state-dependent Hamiltonian. This coupling of the primal and dual variables can have implications on future optimization methods, as we can then attempt various coupling functions, i.e. Hamiltonians, and examine their effectiveness in general optimization techniques.
	
	We also reiterate that our algorithms have been able to perform nonconvex optimization without as many multiple initial guesses as in other algorithm such as coordinate descent. In fact, in all our examples in \cref{sec:examples}, we only used a single initial guess. So we feel a deeper theoretical examination of these algorithms will likely be beneficial to the theory and literature of nonconvex optimization methods.

\newpage
\section{Numerical Results}
\label{sec:examples}

	Here we present numerical examples using our algorithm. The computations were run on a laptop with an Intel i7-4900MQ 2.90Ghz$\times$4 Haswell Core processor, of which only one core processor was utilized for computations. And the computations for the Eikonal equation and the Difference of Norms were computed in C++, while the (unnamed) Isaacs example and the Quadcopter were computed in MATLAB, version R2017b.

	For initializations, we initialized $\tilde{x}^0 = (x_0^0, x_1^0, \ldots, x_N^0)$ to be such that each $x_i^0$ (for $i=0,\ldots,{N-1}$) is a random point close to $x_{\text{target}}$, and we let $x_N^0 \equiv x_{\text{target}}$. In particular, for $\tilde{x}^0$ each coordinate, except for $x_N^0$, was randomly initialized so that $\|\tilde{x}^0 - (x_{\text{target}}, x_{\text{target}}, \ldots, x_{\text{target}})\|_\infty \le 0.1$. We chose $\tilde{p}^0$ to be a random vector close to the origin so that $\|\tilde{p}^0 - 0 \|_\infty \le 0.1$.
	
	We chose $\theta = 1$ in all cases.
	
	The PDHG step-sizes $\sigma$ and $\tau$, and the time-step size $\delta$ all varied for each example.
	
	
	The error tolerance for all optimal control examples were chosen such that primal and dual variables satisfy $\|x^{k+1} - x^k\|_2^2 < 10^{-8}$ and $\|p^{k+1} - p^{k}\|_2^2 < 10^{-8}$. 
	The error tolerance for all the differential games example were also similarly $\|(x^{k+1}, y^{k+1}) - (x^k,y^k)\|_2^2 < 10^{-8}$ and $\|(p^{k+1}, q^{k+1}) - (p^{k},q^k)\|_2^2 < 10^{-8}$, although we had to slightly modify our stopping criteria here. 
	
	For the difference of norms example, if the algorithm reached some max count, then we would examine the value function and stop the algorithm when the value of the value function for consecutive iterations reached a difference below some tolerance.
		
	For the (unnamed) Isaacs example with fully convex initial conditions, we chose the error to be such that the relative error of the value function of consecutive iterations was less than $10^{-6}$, i.e. $\l\|\frac{\text{fval}^{k+1} - \text{fval}^k}{\min(\|\text{fval}^{k+1}\|,1)}\r\| < 10^{-6}$. This example turned out to be the harshest on our algorithm.
	
	In all cases, when we derive the Hamiltonian, we are starting from an optimal control with a terminal condition and solving ``backwards in time" (see \cref{subsec:hjeoc}) as this naturally gives us an initial-valued Hamilton-Jacobi PDE.

\subsection{State-and-Time-Dependent Eikonal Equation (Optimal Control)}
\label{subsec:sateik}

\subsubsection{Brief background on Eikonal equations}
The state-dependent Eikonal equations are HJE with Hamiltonians,
		\begin{align*}
		H^+(x,p) = c(x) \|p\|, \quad H^- = -c(x)\|p\|.
		\end{align*}
	where $c(x) > 0$ and $\|\cdot\|$ is any norm; in our examples we take the Euclidean norm. We also test a state-and-time-dependent equation of Eikonal type, where $H(x,p,t) = c(x-t)\|p\|$. 
	
	They arise from optimal control problems that have dynamics
		\begin{align*}
		f(x,u) = c(x)u, \quad \text{with $\|u\| \leq 1$ and $c(x) \neq 0$ for all $x$}
		\end{align*}
	and with cost-functional
		\begin{align*}
		J[u] = g(x(0)) + \int_0^t \mathcal{I}_{\leq 1}(u(s))\,ds
		\end{align*}
	where $\mathcal{I}_{\leq 1}(\cdot)$ is the indicator function of the unit ball in $\mbb{R}^n$, i.e. 0 for all points within the unit ball including the boundary, and $+\infty$ for all points outside.
	
	This is a \emph{nonlinear} optimal control example due to the presence of $c(x)$. Also, our algorithm is performing nonconvex optimization (due to the presence of $c(x)$, but also in our negative Eikonal equation example where we are performing minimization with a $-g(x)$ term where $g$ is quadratic).
	
	The Eikonal equation features heavily in the level set method \cite{OSHER198812, OshFed03}, which has made wide-ranging contributions in many fields.
	
	
	Note the optimization in solving negative Eikonal equation can be obtained by examining the positive Eikonal equation. This is actually a general phenomenon of Hamiltonians that obey $H(x,-p,t) = H(x,p,t)$. This is because if $\phi$ solves a HJE with initial data $g$ and Hamiltonian such that $H(x,-p,t) = H(x,p,t)$, then examining $-\phi$,
		\begin{align*}
		\l\{\begin{array}{rl}
		(-\phi)_t + H(x, \nabla(-\phi), t) &= 0 \\
		(-\phi)(x,0) &= -g(x)
		\end{array}\r\}
		\quad \Leftrightarrow \quad  
		\l\{\begin{array}{rl}
		\phi_t - H(x,\nabla \phi, t) &= 0 \\
		\phi(x,0) &= g(x) \\
		\end{array}\r\}
		\end{align*}
	so we see $-\phi$ solves the positive Eikonal equation with initial data $-g$ if and only if $\phi$ solve the negative Eikonal equation with initial data $g$. And note that both are viscosity solutions as this computation still holds when we compute the viscous version of the HJE \cite{evans10} (Section 10.1).
	
\subsubsection{Implementation details}
\label{subsubsec:eikimpdet}
	
	Here we take
		\begin{align*}
		c(x) = 1 + 3\text{exp}(-4\|x-(1,1,0,\ldots,0)\|^2_2),
		\end{align*}
	which is a positive bump function. The initial condition of our HJE PDE is
		\begin{align*}
		g(x) = -1/2 + (1/2)\l<A^{-1}x,x\r>, \quad A = \diag(2.5^2, 1, 0.5^2, \ldots, 0.5^2).
		\end{align*}
	
	We used the Lax version of our algorithm, Algorithm \ref{alg:algohjeoclax}, for all cases in this section. For the $\tilde{x}^{k+1}$-update, we took one step of gradient descent. And for the $\tilde{p}^{k+1}$ update, we took the proximal of the $\ell_2$-norm (a.k.a. the $shrink_2$ operator); in general the $shrink$ operator can be defined for any positively homogenous of degree 1 convex function $\phi$. For the $\ell$-1 norm, we have the $shrink_1$ operator (in $\mbb{R}^n$) can be computed coordinate-wise and the $i$-th coordinate satisfies, 
		\begin{align*}
		(\text{shrink}_1(x, \lambda)_i = \l\{ \begin{array}{ll}
		x_i - \lambda & \txt{if } x_i > \lambda \\
		0 & \text{if } |x_i| \le \lambda \\
		x_i + \lambda & \text{if } x_i < -\lambda
		\end{array}\r.
		\end{align*}
	and  the $shrink_2$ operator also can be computed coordinate-wise and the $i$-th coordinate satisfies
		\begin{align*}
		(\text{shrink}_2(x,\lambda))_i = \l\{\begin{array}{ll}
		\frac{x}{\|x\|_2} \max(\|x\|_2 - \lambda, 0) & \text{if } x \neq 0 \\
		0 & \text{if } x = 0.
		\end{array}\r.
		\end{align*}
		
	For the negative Eikonal equation, since in our implementation we computed a positive Eikonal equation with initial data $-g$ and then took the negative, we were able to compute the proximal of the concave quadratic $-g$. We call this taking the $stretch$ operator of $g$ (see \cite{Chow2016}, Section 4.2.2).
		
	For these Eikonal equations, we chose a time step-size of $\delta = 0.02$, and we computed in a $[-3,3]^2$ grid with a mesh size of $0.1$ on each axis. 
	
\begin{figure}[]
		\centering
	    \label{fig:poseik}\includegraphics[scale=0.31]{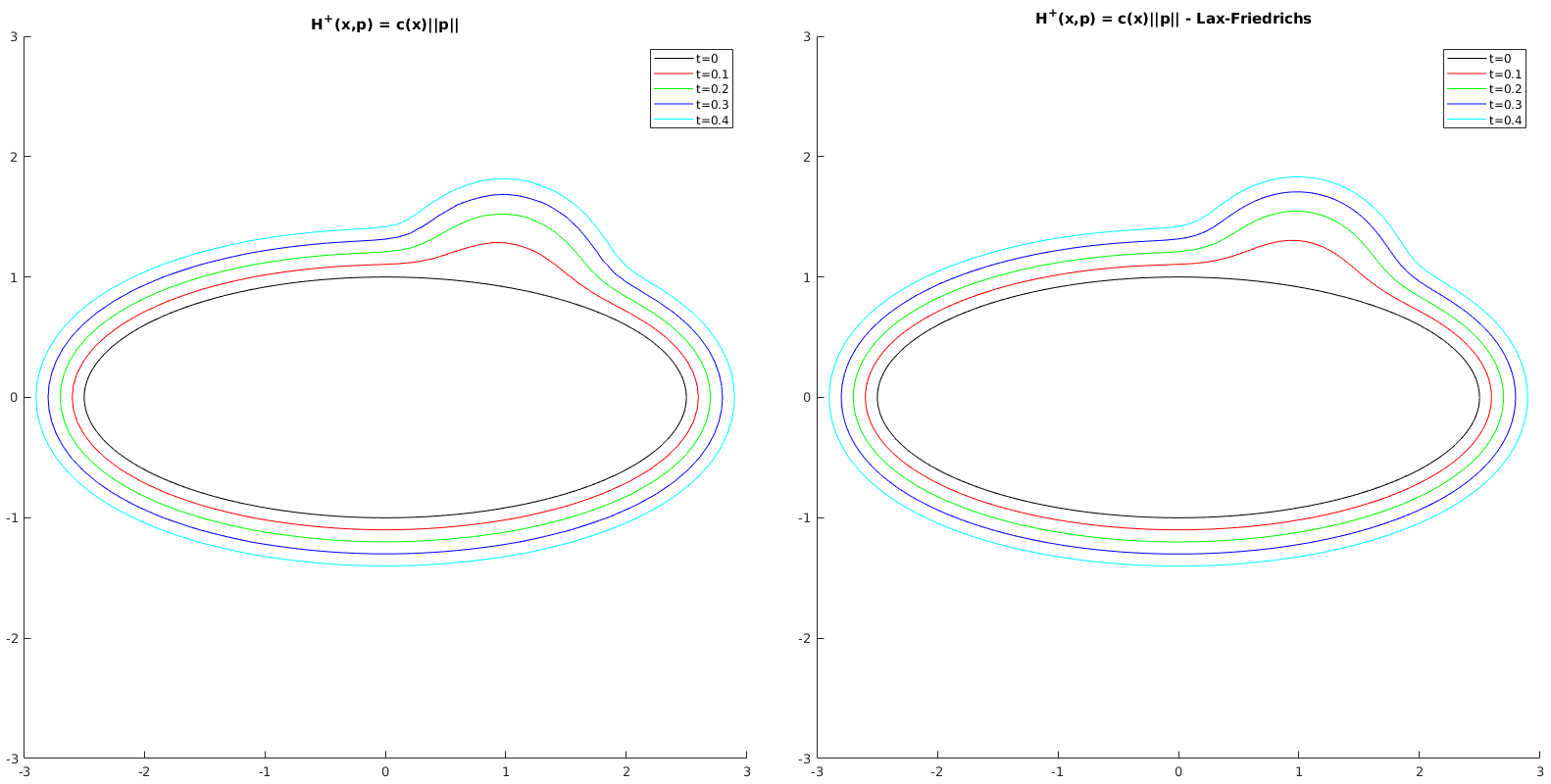}
	    \caption{Eikonal equation with $H^+(x,p) = c(x)\|p\|_2$, in two spatial dimensions. This plot shows the zero level sets of the HJE solution for $t=0.1, 0.2, 0.3, 0.4$. We observe that the zero level sets move outward as time increases. The left figure is computed using our new algorithm, while the right figure is computed using the conventional Lax-Friedrichs method.}
		\end{figure}

	For the positive Eikonal equation (\cref{fig:poseik}), we chose a PDHG step-size that depended on the norm of $\nabla c(x)$. If $\|\nabla c(x)\|_2 > 0.001$, then we took $\sigma = 50$ and or else we took $\sigma = 0.5$. And we always took $\tau = 0.25/\sigma$ (the $0.25$ comes from $\|D\|_2 = 2-\epsilon$ for some small $\epsilon>0$, and requiring $\sigma \tau < 1/\|D\|^2$). To compute the times $t=0.1$, $0.2$, $0.3$, and $0.4$, the computation time averaged to $4.39\times 10^{-4}$  seconds per point in C++.
	
\begin{figure}[]
	    \centering
	    \label{fig:negeik}\includegraphics[scale=0.3]{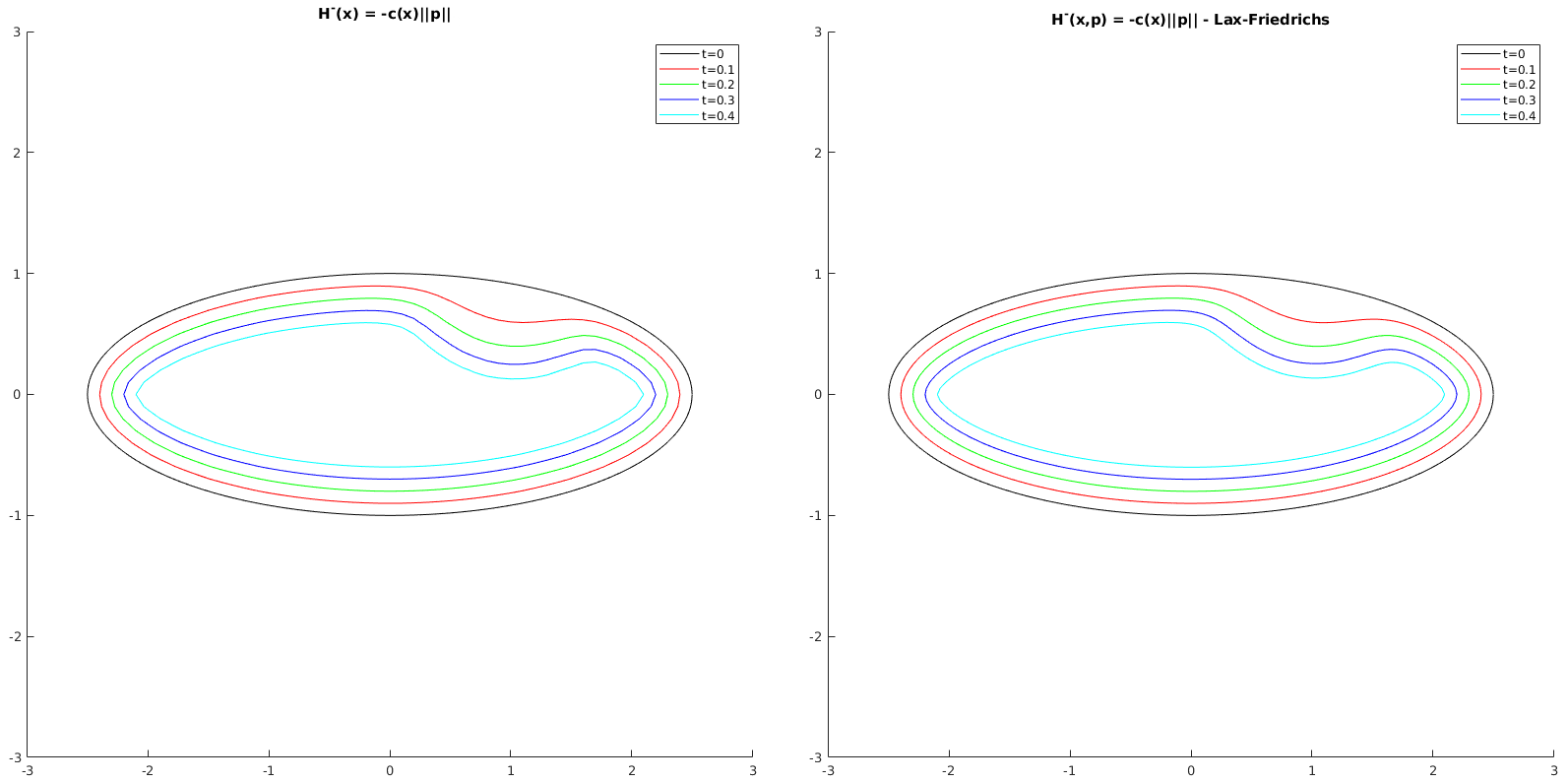}
	    \caption{Eikonal equation with $H^-(x,p) = -c(x)\|p\|_2$, in two spatial dimensions. This plot shows the zero level sets of the HJE solution for $t=0.1, 0.2, 0.3, 0.4$. We observe that the zero level sets move inward as time increases. Left is computed with our new algorithm, while the right is computed using the conventional Lax-Friedrichs method.}
		\end{figure}

	For the negative Eikonal equation (\cref{fig:negeik}), we chose a PDHG step-size of $\sigma = 100$ and $\tau = 0.25/\sigma$. We picked $\theta = 1$ for all cases. For $t=0.1$, $0.2$, $0.3$, $0.4$, and $0.5$, the computation time averaged to 0.0024 seconds per point in C++. We see that for the $t=0.5$ curve, there are kinks, which may be a result of the splitting finding sub/super solutions, rather than viscosity solutions \cite{evans10} (Section 10.1). In this case, multiple initial guesses will alleviate this.
	
\begin{figure}[]
	    \centering
	    \label{fig:negeik10}\includegraphics[scale=0.3]{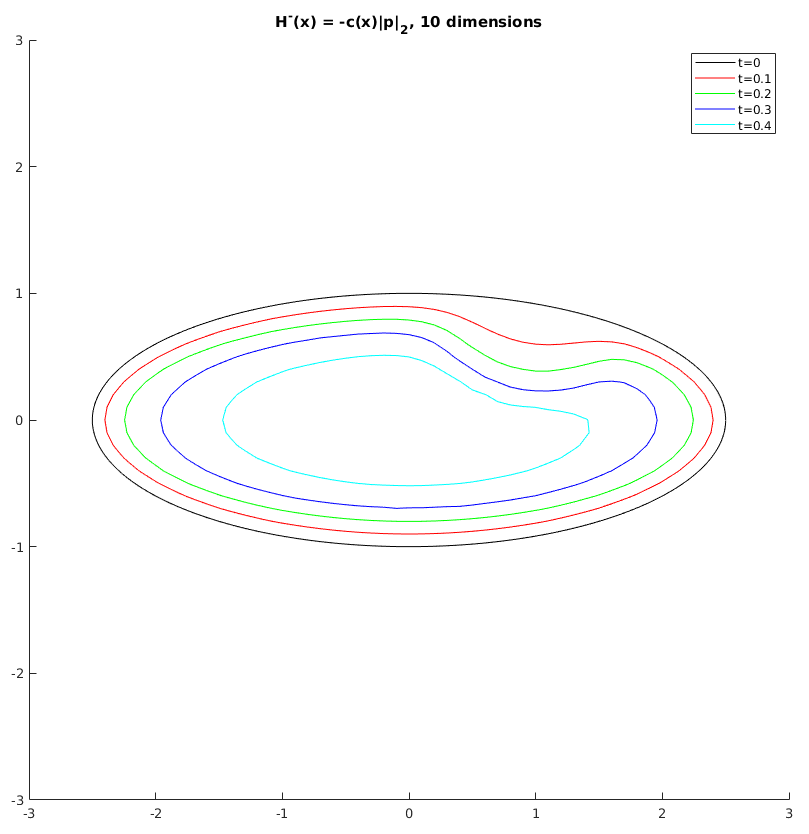}
	    \caption{Eikonal equation with $H^-(x,p) = -c(x)\|p\|_2$, in ten spatial dimensions. This plot shows the zero level sets of the HJE solution for $t=0.1, 0.2, 0.3, 0.4$. We observe that the zero level sets move inward as time increases.}
		\end{figure}

	For the negative eikonal equation in 10 dimensions (\cref{fig:negeik10}), we computed a 2-dimensional slice in $[-3,3]\times [-3,3] \times \{0\} \times \cdots \times \{0\}$. The time step-size was again $\delta = 0.02$ and we chose $\sigma = 100$, and $\tau = 0.25/\sigma$. For $t=0.1$, $0.2$, $0.3$, and $0.4$ ($0.5$ had no level sets) the computation time averaged to 0.004 seconds per point in C++.
	
\begin{figure}[]
	    \centering
	    \label{fig:poseiktime}\includegraphics[scale=0.3]{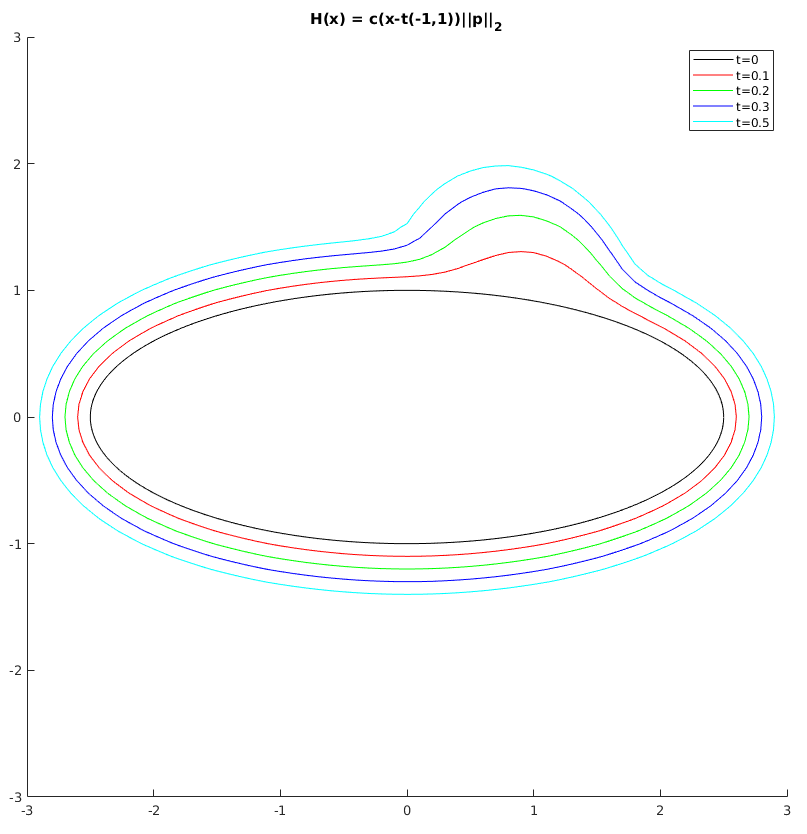}
	    \caption{Eikonal equation with $H(x,p) = c(x-t(-1,1))\|p\|_2$. This plot shows the zero level sets of the HJE solution for $t=0.1, 0.2, 0.3, 0.4$. We observe there is a similarity to the positive eikonal case, but the ``bump" is more sheared to the left.}
		\end{figure}

	We also computed a state-and-time-dependent Eikonal equation (\cref{fig:poseiktime}) where $H(x,p,t) = c(x-t(-1,1))\|p\|_2$. Here in our specific example, $c(x-t(-1,1))$ represents a bump function moving diagonally in the $(-1,1)$ direction as $t$ increases. The time step-size were the same as in the positive eikonal case which is reasonably expected because we took the positive eikonal case and modified it. The computational time for $t=0.1$, $0.2$, $0.3$, and $0.4$ averaged to $5.012 \times 10^{-4}$ seconds per point in C++.
		
	In these examples, we achieved a speedup of about ten times over coordinate descent, and only one initial guess was used. In low dimensions, this problem can be solved with SQP (Sequential Quadratic Programming) on the value function, or using Lax-Friedrichs. We use these methods to check our accuracy and they agree to within $10^{-4}$ for each point $(x,t)$ when using SQP.

	We observe that the two different eikonal equations required vastly different step-sizes and it is worth examining how to choose step-sizes in a future work. We speculate the step-sizes may act as CFL conditions. This is a further point of study.
	
	
	Comparing coordinate descent to our algorithm in MATLAB, we achieve about an 8-10 times speed-up. 
	
	The figures \cref{fig:poseik}, \cref{fig:negeik}, and \cref{fig:negeik10}, and \cref{fig:poseiktime} show the zero level sets of the HJE solution.
	
\subsubsection{Dimensional scaling of the negative Eikonal equation}
\label{subsubsec:dimscale}
	
	Here we examine how Algorithm 1 scales with dimension. We compute the negative Eikonal equation with the same speed $c$ and initial data $g$ as above, and with $\delta = 0.02$ and $\sigma = 100$ and $\tau = 0.25/\sigma$. We computed in a 2-dimensional slice $[-3,3]^2 \times \{0\}^{d-2}$, and we computed from $d=10$ to $d=2000$ dimensions. We performed our analysis at time $t=0.2$.
	
	We chose this particular example as this is a nonlinear optimal control problem that requires us to optimize a nonconvex problem.
	
\begin{figure}[]
	    \centering
	    \label{fig:dimscale}\includegraphics[scale=0.5]{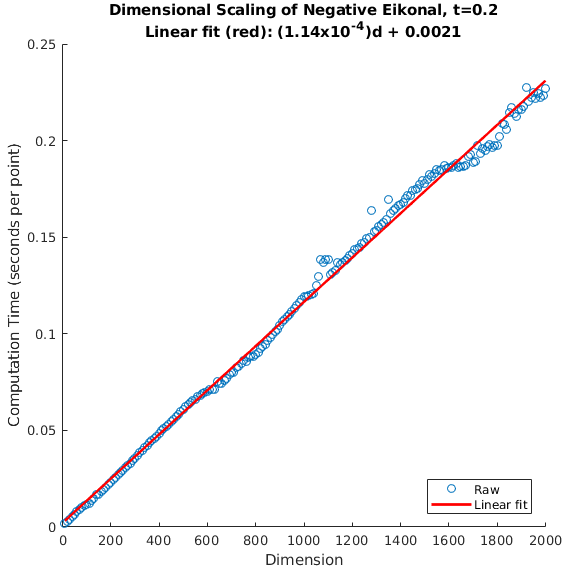}
	    \caption{How our algorithm scales with dimension for the negative Eikonal equation at time $t=0.2$. This is a nonlinear optimal control problem, and the optimization requires us to perform nonconvex optimization. The plot shows a linear fit.}
		\end{figure}

	We used least-squares to fit both a linear function as well as a quadratic function. The coefficients were
		\begin{align*}
		\text{lin}(d) = (1.14\times 10^{-4}) d + 0.0021, \quad \text{quad}(d) = (-5.99 \times 10^{-9}) d^2 +  (1.27 \times 10^{-4}) d - 0.00195
		\end{align*}	
	As we can see from the equations of the fit, and from \cref{fig:dimscale}, the quadratic fit has an extremely small quadratic coefficient. \cref{fig:dimscale} shows the plot with the linear fit. This computation was done in C++.
	
	We predict that for general problems, the scaling will be polynomial in time.
	
\subsection{Difference of norms (Differential Games)}

\subsubsection{From differential games to HJE for the difference of norms}

	The state-dependent HJE for the difference of norms case arises from the following differential games problem: Given $x\in \mbb{R}^{d_1}$ and $y\in \mbb{R}^{d_2}$, and some $t>0$, we have the following the dynamics
		\begin{align*}
		\l\{\begin{array}{l}
		\begin{pmatrix} \dot{\tbf{x}}(s) \\ \dot{\tbf{y}(s)} \end{pmatrix} = \begin{pmatrix} c_1(\tbf{x}(s), \tbf{y}(s)) \, \bsym{\alpha}(s) \\ c_2(\tbf{x}(s),\tbf{y}(s)) \, \bsym{\beta}(s)\end{pmatrix}, \quad 0 \le s \le t\\ 
		\begin{pmatrix} \tbf{x}(t) \\ \tbf{y}(t) \end{pmatrix} = \begin{pmatrix} x \\ y \end{pmatrix} \\ \\
		\bsym{\alpha}(s)\| \le 1, \|\bsym{\beta}(s)\| \le 1, \text{ for all $s$},\\ \\
		\text{and } c_1, c_2 \text{ are positive functions.}  
		\end{array}\r.
		\end{align*}
	And the cost function is
		\begin{align*}
		J_{x, t}[\bsym{\alpha}, \bsym{\beta}] = g(\tbf{x}(0), \tbf{y}(0)) + \int_0^t \mathcal{I}_{\le 1}(\bsym{\alpha}(s)) - \mathcal{I}_{\le 1}(\bsym{\beta}(s)) \, ds
		\end{align*}
	where $\mathcal{I}_{\le 1}(\cdot)$ is the indicator function the unit-ball, i.e. it equals $0$ for points inside and on the unit-ball, and $+\infty$ for points outside. Our value function is then 
		\begin{align*}		
		\phi(x,y,t) = \inf_{\alpha, \|\alpha\|\le 1}\sup_{\beta, \|\beta\| \le 1} J_{x,t}[\bsym{\alpha}, \bsym{\beta}].
		\end{align*}
	Then our Hamiltonian becomes,
		\begin{align*}
		H(x,y,p,q) &= \max_{\alpha} \min_{\beta} \l\{ \l<\begin{pmatrix}c_1(x,y) \alpha \\ c_2(x,y)\beta\end{pmatrix}, \begin{pmatrix}p \\ q\end{pmatrix} \r> - \l(I_{\le 1}(\alpha) - I_{\le 1}(\beta) \r) \r\} \\
		&= \max_{\alpha} \l\{ \l<c_1(x,y)\alpha, p_1\r> - I_{\le 1}(\alpha)\r\} + \min_{\beta}  \l\{ \l<c_2(x,y)\beta, q\r> +  I_{\le 1}(\beta) \r\} \\
		&= c_1(x,y)\|p\|_2 - c_2(x,y) \|q\|_2.
		\end{align*}
		
	In this case, we have nonlinear dynamics, aswell as nonconvex Hamiltonian.
	
	\subsubsection{Implementation details for the difference of norms}
	
	Here we take,
		\begin{align*}
		c_1(x) = 1 + 3\text{exp}(-4\|x-(1,1,0,\ldots,0)\|^2_2), \qquad c_2(x) = c_1(-x),
		\end{align*}
	and the initial condition is
		\begin{align*}
		g(x) = -1/2 + (1/2)\l<A^{-1}x,x\r>, \quad A = \diag(2.5^2, 1, 0.5^2, \ldots, 0.5^2).
		\end{align*}
	which is the same initial condition as in our Eikonal equation example above \cref{subsubsec:eikimpdet}.
	
	For the $2$-dimensional case, we used the Hamiltonian,
		\begin{align*}
		H(x_1, x_2, p_1, p_2) = c(x_1, x_2) \|p_1\|_2 - c(-x_1, -x_2) \|p_2\|_2
		\end{align*}
	and for the $7$-dimensional case, we used,
		\begin{align*}
		H(x_1, x_2,\ldots, x_7, p_1, p_{(2,\ldots,7)}) = c(x_1, \ldots, x_7) \|p_1\|_2 - c(-x_1, \ldots, -x_7) \|p_{(2,\ldots,7)}\|_2.
		\end{align*}
	where $p_{(2,\ldots,7)} = (p_2, p_3, \ldots, p_7)$.
	
	We compute our solutions in $2$ and $7$ dimensions. We compute the $2$-dimensional case in a $[-3,3]^2$ grid, and we compute the $7$-dimensional case in the two dimensional slice $[-3,3]^2 \times \{0\} \times \cdots \times \{0\}$. And we used \cref{alg:algohjedglax}.
	
	For the $\tilde{p}^{k+1}$-update, we used the proximal of the $\ell_2$-norm (a.k.a. the $shrink_2$ operator), and for the $\tilde{x}^{k+1}$-update, we used one step of gradient descent.
	
	We took the time-step as $\delta = 0.02$, and we took the PDHG steps $\sigma = 50$, and $\tau = 0.25/\sigma$. The $0.25$ comes from the fact that the PDHG algorithm requires $\sigma \tau \|D\|_2^2 < 1$, and $\|D\|_2 = 2 - \epsilon$, for some small $\epsilon > 0$.

	The computation was done with a mesh size of $1/12\approx 0.08333$ in each axis. For the 2-dimensional case, the computation averaged out to 0.0135 seconds per point in C++, and for the 7-dimensional case the computation averaged out to 0.01587 seconds per point in C++. If we compared the algorithms in MATLAB on the same computer, we achieved a 10-20 times speed-up compared to coordinate descent.
	
	We note that at certain points, the trajectories would oscillate a little for larger times which may be due to the non-convexity and the non-unique optimal trajectories. So when a maximum count was reached, we would raise $\sigma$ by 20, and we would also readjust $\tau = 0.25/\sigma$. We do not recommend choosing a high $\sigma$ at all points, or else the algorithm would result in incorrect solutions. If the convergence was not fast enough, after some maximum count, we would switch our convergence criteria to the value function, as the error (between consecutive iterations) seemed to converge to zero. The procedure of raising the sigma is reminiscent of CFL conditions for finite-difference schemes, and we are examining how best to choose the PDHG step-sizes $\sigma$ and $\tau$. The best $\sigma$ and $\tau$ to choose seem to be dependent on the point at which we are computing.
	
	
	
		\begin{figure}[]
		\centering
	    \label{fig:2dimdiffgames}\includegraphics[scale=0.31]{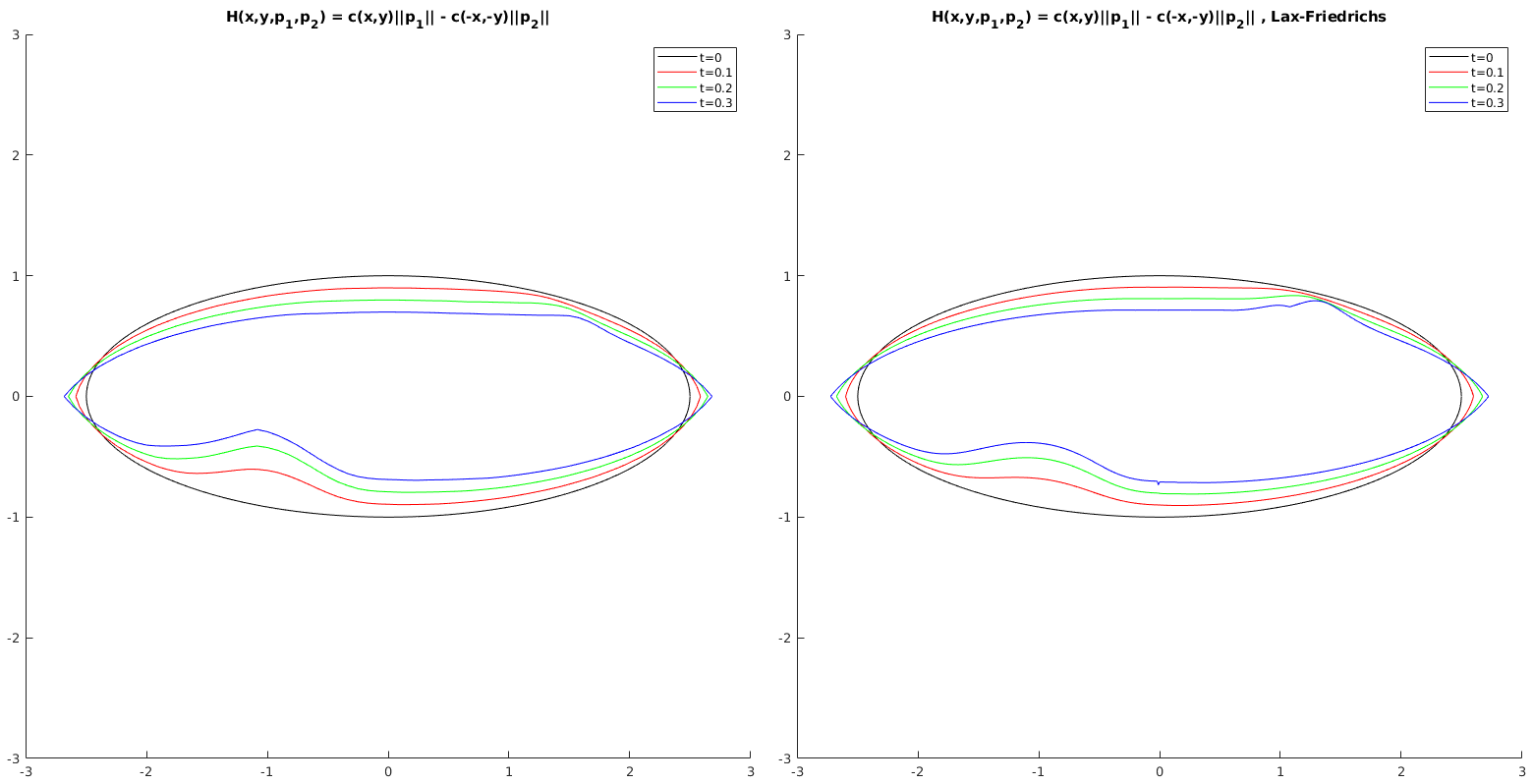}
	    \caption{The difference of norms HJE in two spatial dimensions. This plot shows the zero level sets of the HJE solution for $t=0.1, 0.2, 0.3$. We observe that the zero level sets move inward as time increases. Left is computed with our new algorithm, while the right is computed using the conventional Lax-Friedrichs method. Note there is an anomaly at the top-right of Lax-Friedrich computation. And there is also more of a corner in the bottom-left of the solution computed by the new method. This may be a result of the true solution, and which does not appear in the Lax-Friedrichs solution as it tends to smooth out solutions.}
		\end{figure}
		
		\begin{figure}[]
		\centering
	    \label{fig:7dimdiffgames}\includegraphics[scale=0.31]{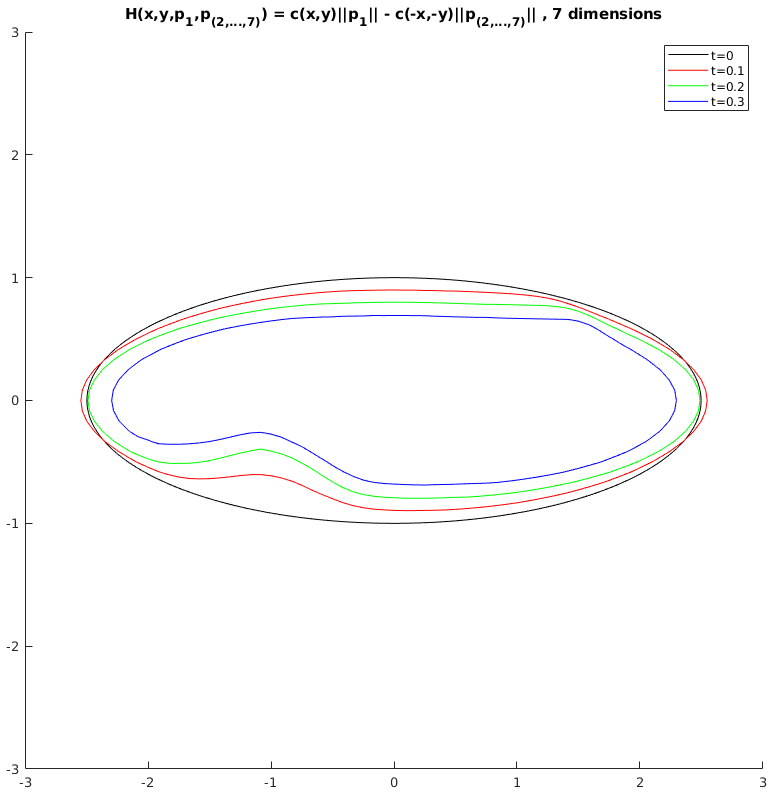}
	    \caption{The difference of norms HJE in seven spatial dimensions. This plot shows the zero level sets of the HJE solution for $t=0.1, 0.2, 0.3$.}
		\end{figure}

\subsection{An (unnamed) example from Isaacs (Differential Games)}
\label{subsec:unnamedIsaacs}

\subsubsection{From differential games to HJE for an unnamed example from Isaacs}

	We modify an example from \cite{evansdgslides}, to obtain the following example of a differential game. The dynamics are as follows:
		\begin{align*}
		\l\{\begin{array}{rl}
		\dot{x}(s) &= 2 \beta + \sin(\alpha) \\ 
		\dot{y}(x) &= -c(x,y) + \cos(\alpha)
		\end{array}\r.
		\end{align*}
	where $0 \le \alpha \le 2\pi$ and $-1 \le \beta \le 1$. These dynamics are nonlinear. We take the cost-functional as,
		\begin{align*}
		J[\alpha,\beta] = g(x(0), y(0)) + \int_{0}^t 1\, ds
		\end{align*}
	and the value function seeks to maximize with respect to $\alpha \in [0,2\pi]$, and minimize with respect to $\beta \in [-1,1]$. Then our Hamiltonian is,
		\begin{align*}
		H(x,y,p,q) &= \min_{\alpha, \alpha \in [0,2\pi]} \max_{\beta, \beta \in [-1,1]} \l\{ \l< \begin{pmatrix}2 \beta + \sin(\alpha) \\ -c(x,y) + \cos(\alpha)\end{pmatrix}, \begin{pmatrix} p \\ q \end{pmatrix}\r> - 1\r\} \\
		&= \min_{\alpha, \alpha \in [0,2\pi]} \max_{\beta, \beta \in [-1,1]} \l\{ 2\beta p - c(x,y) q + p\sin(\alpha) + q\cos(\alpha) - 1\r\} \\
		&= -c(x,y)q + 2|p| - \sqrt{p^2 + q^2} - 1.
		\end{align*}			
	This Hamiltonian is nonconvex. And the dynamics are nonlinear.
	
\subsubsection{Implementation details for the (unnamed) example from Isaacs with fully convex initial conditions}

	We take
		\begin{align}\label{eq:cisaacsfullyconvex}
		c(x) = 2( 1 + 3\text{exp}(-4\|x-(1,1,0,\ldots,0)\|^2_2) ),
		\end{align}
	which is a positive bump function. The initial condition of our HJE PDE is
		\begin{align*}
		g(x) = -1/2 + (1/2)\l<A^{-1}x,x\r>, \quad A = \diag(2.5^2, 1, 0.5^2, \ldots, 0.5^2).
		\end{align*}
	
	This example is perhaps the harshest on our algorithm and turns out to be slower than coordinate descent, but in many ways this is not surprising. This is because this problem is highly nonconvex and the $g(x,y)$ that we use is a convex function -- ideally we would like it to be a convex-concave function which would be suitable for saddle-point problems. Not only that, but our Hamiltonian is not bounded below with respect to $q$, and the Hopf formula requires this assumption.
	
	Nevertheless, we show this example in order to advertise the generality of our algorithm. It might actually be surprising that our algorithm gives a solution that looks like the Lax-Friedrichs solution at all. We also not that we only used one initial guess, and in our experiments, using around $5$ initial guesses smooths our the solution.
	
	We use \cref{alg:algohjedglax}, but modified so we can utilize the convex portion of $g$ (see the last paragraph of \cref{subsec:splitforhjedg}). We compute our solutions in a $2$-dimensional $[-3,3]^2$ grid. The $\tilde{p}^{k+1}$-update utilizes a combination of gradient descent for the $q$, and the $shrink_2$-operator the $p$. The $\tilde{x}^{k+1}$-update uses one step of gradient descent.
	
	We took the time-step as $\delta = 0.005$, and we took the PDHG steps as $\sigma = 20$, and $\tau = 0.25/\sigma$, where the $0.25$ comes from the PDHG condition that $\sigma \tau \|D\|_2^2 < 1$, and $\|D\|_2 = 2 - \epsilon$ for some small $\epsilon > 0$.
	
	As in the difference of norms example, we did have some points that were slower to converge. We alleviated this problem by rerunning the algorithm at the same point (without change $\sigma$ nor $\tau$) if we reached some maximum count, and sometimes took the solution if the maximum count was reached anyway. 
	
	The computational time on a $[-3,3]^2$ grid, of mesh size $1/12\approx 0.0833$ for each axis, averaged out to 0.412 seconds per point in MATLAB.	
	
	In this example, using \cref{alg:algohjedglax} was slower than coordinate descent where it averaged to 0.133 seconds per point, and so we recommend using coordinate descent in this case.
	
\begin{figure}[]
		\centering
	    \label{fig:isaacsplot}\includegraphics[scale=0.31]{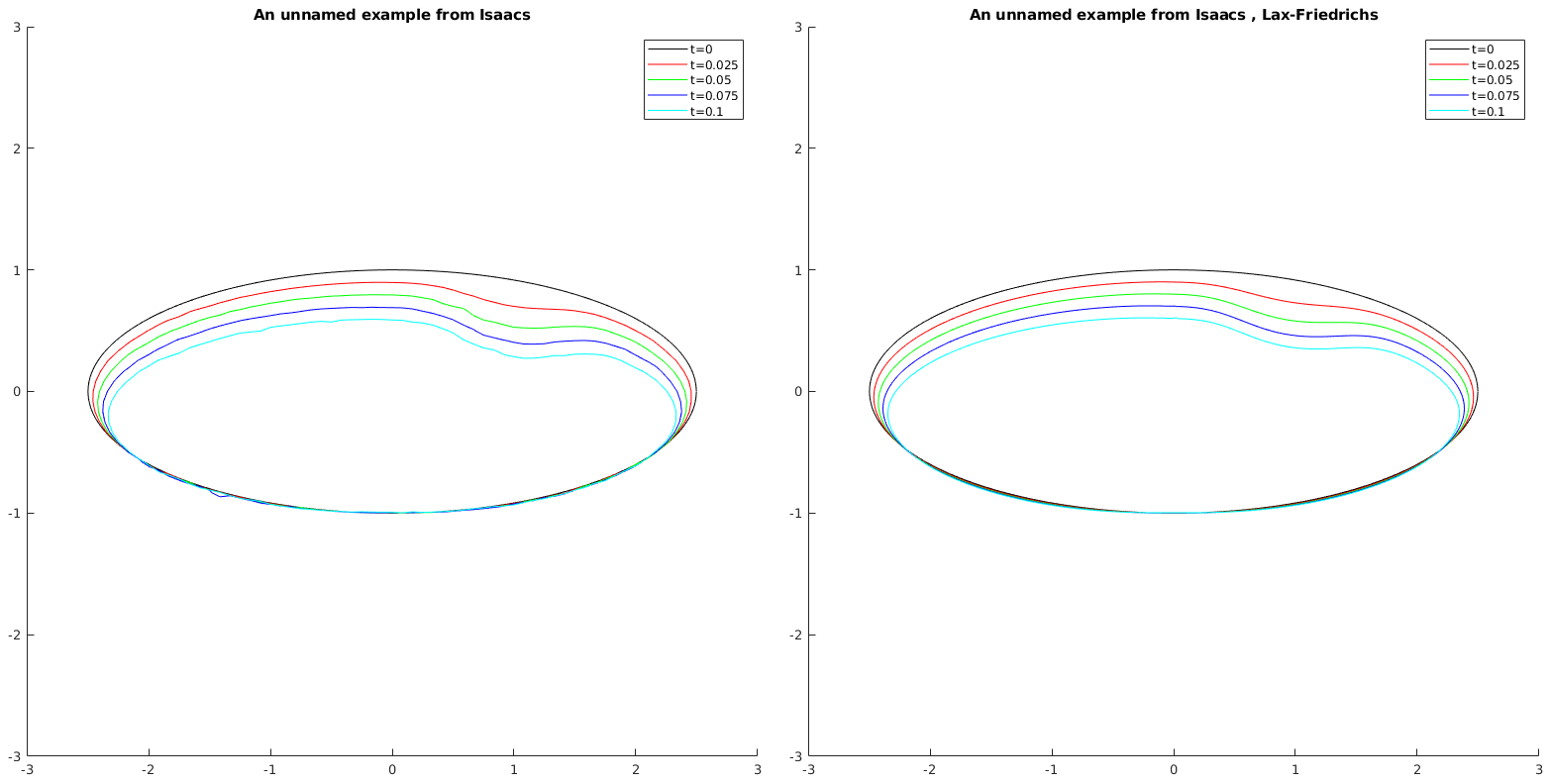}
	    \caption{The zero level-sets for the HJE from an unnamed example of Isaacs. The times we computed were $t=0.025$, $0.05$, $0.075$, and $0.1$. The left figure is the result of our algorithm, while the right figure is the result of Lax-Friedrichs. Here we see this example is our harshest critic. But this is not surprising because the initial condition is a fully convex function, whereas we'd rather have it be convex-concave. And also the Hamiltonian is not bounded below with respect to $q$ which as an assumption of the Hopf formula. Neverthless our algorithm is still able to achieve a result similar to Lax-Friedrichs, which might actually be the surprising part. We also note that we only used one initial guess here, but using multiple initial guess (around 5) smooths out the curves.}
		\end{figure}

	\cref{fig:isaacsplot} gives the result of our algorithm.
	
\subsubsection{Implementation details for the (unnamed) example from Isaacs with convex-concave initial conditions}
\label{subsubsec:isaacscc}

	We take $c(x)$ to be the same as in the fully convex initial conditions (see \cref{eq:cisaacsfullyconvex}). The initial condition of our HJE PDE is
		\begin{align*}
		g(x_1, x_2) = -1/2 + (1/2)( \l((2.5) x_1\r)^2 - \l(x_2\r)^2)
		\end{align*}
	
	Here we have convex-concave initial conditions, and our algorithm works well.
	
	We use \cref{alg:algohjedghopf}, and as in all other examples here, we only have one initial guess. We compute our solutions in a $2$-dimensional $[-3,3]^2$ grid. The $\tilde{p}^{k+1}$-update utilizes a combination of gradient descent for the $q$, and the $shrink_2$-operator the $p$. The $\tilde{x}^{k+1}$-update uses one step of gradient descent.
	
	We took the time-step as $\delta = 0.005$, and we took the PDHG steps as $\sigma = 2$ for $t=0.025$, $0.05$, $0.075$, and $\sigma = 10$ for $t=0.1$. We always chose $\tau = 0.25/\sigma$, where the $0.25$ comes from the PDHG condition that $\sigma \tau \|D\|_2^2 < 1$, and $\|D\|_2 = 2 - \epsilon$ for some small $\epsilon > 0$.
	
	We computed this example in a $[-3,3]^2$ grid with mesh size $1/12 \approx 0.0833$ for each axis. The computational time averaged to 0.125 seconds per point.
	
	As in the difference of norms example, we did have some points that were slower to converge. We alleviated this problem by rerunning the algorithm at the same point (without change $\sigma$ nor $\tau$) if we reached some maximum count, and sometimes took the solution if the maximum count was reached anyway.
	
\begin{figure}[]
		\centering
	    \label{fig:isaacsconvexconcaveplot}\includegraphics[scale=0.31]{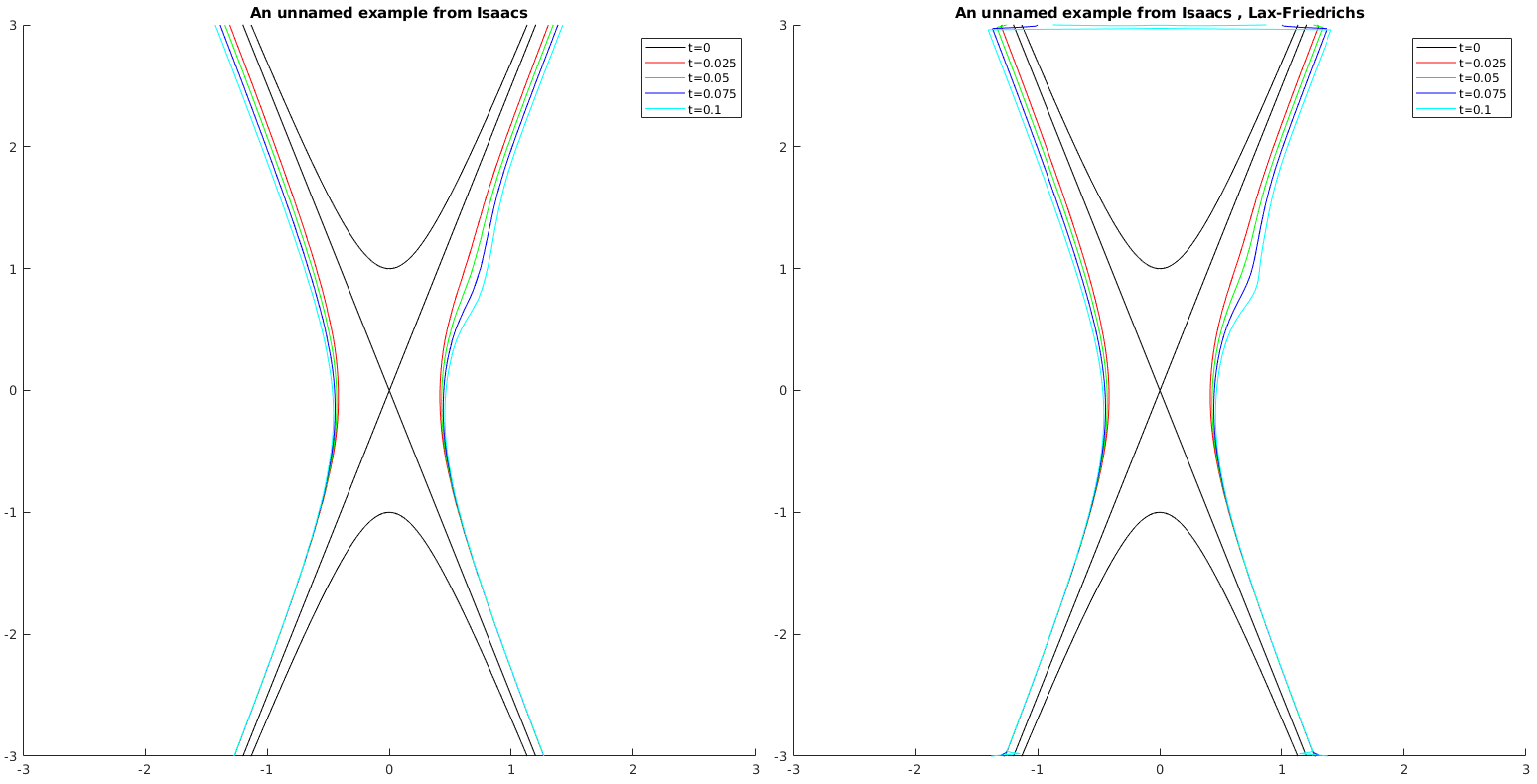}
	    \caption{The zero level-sets for the HJE from an (unnamed) example of Isaacs with convex-concave initial conditions. The times we computed were $t=0.025$, $0.05$, $0.075$, and $0.1$. The left figure is the result of our algorithm, while the right figure is the result of Lax-Friedrichs.}
		\end{figure}

	\cref{fig:isaacsconvexconcaveplot} gives the result of our algorithm. 
	
\subsection{Quadcopter (a.k.a. Quadrotor or Quad rotorcraft) (Optimal Control)}
\label{sec:quadcopter}

\subsubsection{From optimal control to HJE for the quadcopter}

	A quadcopter is a multirotor helicopter that utilizes four rotors to propel itself across space. The dynamics of a quadcopter \cite{GarciaCarrillo2013} are:
		\begin{align*}
		\l\{ \begin{array}{rl}
		\ddot{x} &= \frac{u}{m}\l( \sin(\phi)\sin(\psi) + \cos(\phi)\cos(\psi)\sin(\theta) \r) \\ 
		\ddot{y} &= \frac{u}{m}\l(-\cos(\psi)\sin(\phi) + \cos(\phi)\sin(\theta)\sin(\psi) \r) \\ 
		\ddot{z} &= \frac{u}{m}\cos(\theta)\cos(\phi) - g \\ 
		\ddot{\psi} &= \tilde{\tau}_{\psi} \\ 
		\ddot{\theta} &= \tilde{\tau}_{\theta} \\ 
		\ddot{\phi} &= \tilde{\tau}_{\phi}
		\end{array} \r.
		\end{align*}
	where $(x,y,z)$ is the position of the quadcopter in space, and $(\psi,\theta,\phi)$ is the angular orientation of the quadcopter (a.k.a. Euler angles). The above second-order system turns into the first-order system,
		\begin{align*}
		\l\{ \begin{array}{rl}
		\dot{x}_1 &= x_2 \\ 
		\dot{y}_1 &= y_2 \\ 
		\dot{z}_1 &= z_2 \\ 
		\dot{\psi}_1 &= \psi_2 \\ 
		\dot{\theta}_1 &= \theta_2 \\ 
		\dot{\phi}_1 &= \phi_2 \\ 
		\dot{x}_2 &= \frac{u}{m}\l( \sin(\phi_1)\sin(\psi_1) + \cos(\phi_1)\cos(\psi_1)\sin(\theta_1) \r) \\ 
		\dot{y}_2 &= \frac{u}{m}\l(-\cos(\psi_1)\sin(\phi_1) + \cos(\phi_1)\sin(\theta_1)\sin(\psi_1) \r) \\ 
		\dot{z}_2 &= \frac{u}{m}\cos(\theta_1)\cos(\phi_1) - g \\ 
		\dot{\psi}_2 &= \tilde{\tau}_{\psi} \\ 
		\dot{\theta}_2 &= \tilde{\tau}_{\theta} \\ 
		\dot{\phi}_2 &= \tilde{\tau}_{\phi}
		\end{array} \r.
		\end{align*}
	and so the right-side becomes our $\tbf{f}(\tbf{x},\bsym{\alpha})$. Here the controls are the variable $u$, $\tilde{\tau}_\psi$, $\tilde{\tau}_\theta$, $\tilde{\tau}_{\phi}$. 
	
	This is a 12-dimensional, nonlinear, optimal control problem.
	
	Denoting $\tbf{x} = (x_1, y_1, z_1, \psi_1, \theta_1, \phi_1, x_2, y_2, z_2, \psi_2, \theta_2, \phi_2)$, then our cost-functional is,
		\begin{align}\label{eq:quadcostfunc}
		J[u, \tilde{\tau}_{\psi}, \tilde{\tau}_{\theta}, \tilde{\tau}_{\phi}] = g(\tbf{x}(0)) + \int_0^t 2 + \|(u(s), \tilde{\tau}_{\psi}(s), \tilde{\tau}_{\theta}(s), \tilde{\tau}_{\phi}(s))\|_2^2 \,ds 
		\end{align}
	where this cost functional was chosen to follow \cite{7759553} and \cite{7040310}. Therefore, our Hamiltonian becomes,
		{\small
		\begin{align*}
		H(\tbf{x}, \tbf{p}, t) &= \max_{u, \tilde{\tau}_{\psi}, \tilde{\tau}_\theta, \tilde{\tau}_{\phi} } \l\{ \begin{bmatrix}x_2 \\ y_2 \\ z_2\end{bmatrix} \cdot \begin{bmatrix}p_1 \\ p_2 \\ p_3\end{bmatrix}  + \begin{bmatrix}\psi_2 \\ \theta_2 \\ \phi_2\end{bmatrix} \cdot \begin{bmatrix}p_4 \\ p_5 \\ p_6\end{bmatrix} \r. \\ 
		& \qquad \l. + \frac{u}{m}\begin{bmatrix}\sin(\phi_1)\sin(\psi_1) + \cos(\phi_1)\cos(\psi_1)\sin(\theta_1) \\ -\cos(\psi_1)\sin(\phi_1) + \cos(\phi_1)\sin(\theta_1)\sin(\psi_1) \\ \cos(\theta_1)\cos(\phi_1)\end{bmatrix} \cdot \begin{bmatrix}p_7 \\ p_8 \\ p_9\end{bmatrix} - p_9g + \begin{bmatrix}\tilde{\tau}_{\psi} \\ \tilde{\tau}_{\theta} \\ \tilde{\tau}_{\phi}\end{bmatrix} \cdot \begin{bmatrix}p_{10} \\ p_{11} \\ p_{12}\end{bmatrix}\r. \\
		& \qquad \l. - 2 - ||u||^2 - ||\tilde{\tau}_{\psi}||^2 - ||\tilde{\tau}_{\theta}||^2 - ||\tilde{\tau}_{\phi}||^2\r\} \\ \\
		&=  \begin{bmatrix}x_2 \\ y_2 \\ z_2\end{bmatrix} \cdot \begin{bmatrix}p_1 \\ p_2 \\ p_3\end{bmatrix}  + \begin{bmatrix}\psi_2 \\ \theta_2 \\ \phi_2\end{bmatrix} \cdot \begin{bmatrix}p_4 \\ p_5 \\ p_6\end{bmatrix} \\
		& + \frac{1}{4 m} \l\|\begin{bmatrix}\sin(\phi_1)\sin(\psi_1) + \cos(\phi_1)\cos(\psi_1)\sin(\theta_1) \\ -\cos(\psi_1)\sin(\phi_1) + \cos(\phi_1)\sin(\theta_1)\sin(\psi_1) \\ \cos(\theta_1)\cos(\phi_1)\end{bmatrix} \cdot \begin{bmatrix}p_7 \\ p_8 \\ p_9\end{bmatrix} \r\|^2 \\
		& - p_9g + \frac{1}{4}\|p_{10}\|^2 + \frac{1}{4}\|p_{11}\|^2 + \frac{1}{4}\|p_{12}\|^2 - 2
		\end{align*}
	}
	
	\subsubsection{Implementation details for the quadcopter}
	
	Here we have,
		\begin{align*}
		g(x) = -1/2 + (1/2)\l<A^{-1}x,x\r>, \quad A = \diag(0.2, 1, 1, \ldots, 1).
		\end{align*}
	In this case, we use the algorithm based on the generalized Hopf formula, \cref{alg:algohjeochopf}. 
	
	We compute our solutions in a two dimensional slice of $\mbb{R}^{12}$:
		\begin{align*}
		[-1,1] \times \{0\} \times \{0\} \times [-1,1] \times \{0\} \times \{0\} \times \{0\}\times \{0\} \times \{0\} \times \{0\} \times \{0\} \times \{0\},
		\end{align*}
	i.e. we vary the $x$-coordinate, as well as the $x$-velocity-coordinate. Recall the order of the coordinates are: $\tbf{x} = (x_1, y_1, z_1, \psi_1, \theta_1, \phi_1, x_2, y_2, z_2, \psi_2, \theta_2, \phi_2)$.

	For both the $\tilde{p}^{k+1}$-update and the $\tilde{x}^{k+1}$-update, we used gradient descent, except for the update involving $g^*(p_1)$, where we did a proximal-gradient step, i.e. we performed a gradient descent, ignoring $g^*$, and then we fed the result into $\text{prox}_{\sigma (g^*)}(\cdot)$.
	
	In this example, we chose as time-step size $\delta = 0.005$, and we chose times $t=0.025, 0.05, 0.075$. For the PDHG step-sizes, we chose $\sigma = 5$, and $\tau = 0.25/\tau$, where as stated above, the $0.25$ comes from the PDHG requirement $\sigma \tau \|D\|_2^2 < 1$, and $\|D\|_2 = 2 - \epsilon$, for some small $\epsilon>0$.
	
	The computation was done on a $[-1,1]^2$ grid, with mesh size $0.01$ in each axis. The computational time averaged to 0.0733 seconds per point in MATLAB.
	
	In this case, not only are we able to compute level-sets of the HJE, but we are also able to take advantage of the characteristic curve/optimal trajectory generation that is freely offered by our algorithm.
	
	To generate the curves/trajectories, we took a randomly-generated terminal point, which was exactly
		\begin{align*}
		x_{\txt{target}} = (0.36, -0.62, -0.06, 0.23, 0.85, -0.66, 0.72, -0.45, 0.15, -0.75, 0.04, -0.83)
		\end{align*}				
	and we computed up to $t=6$ seconds. We chose a time-step of $\delta = 0.05$, and we chose $\sigma = 11$ (and $\tau = 0.24/\sigma$ (as opposed to $0.25$ as the latter will not converge). This took about 24s to compute in MATLAB. We verified our result by directly minimizing a discretized version of \eqref{eq:quadcostfunc} (see also \cref{sec:discgenlaxhopf} on how we discretized), which is a direct collocation method \cite{vonStryk1993}. We performed the optimization using a standard MATLAB minimization solver (fmincon with the SQP algorithm), and this agreed with our results. The solver took 133-347s to compute the trajectories, depending on the accuracy criteria, and we note that fmincon converges to our splitting result the longer we let the algorithm run. So in this case, 5-10+ times speedup. Computing trajectories at other points have found around an 8-10+ speedups.
	
\begin{figure}[]
		\centering
	    \label{fig:quadcopterhje}\includegraphics[scale=0.31]{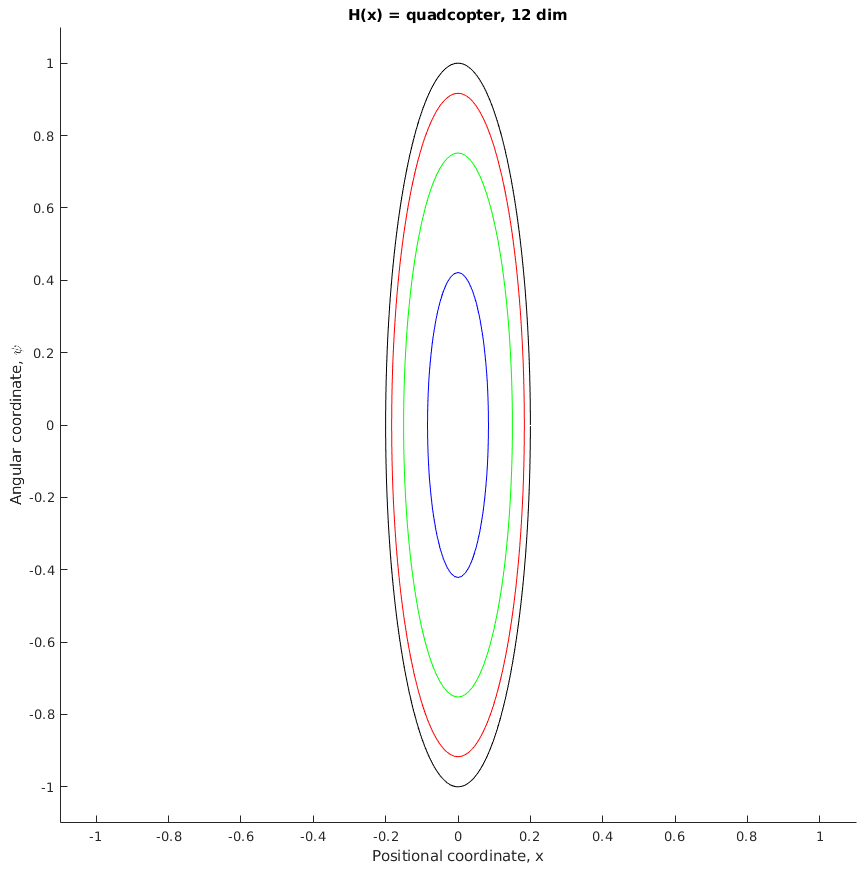}
	    \caption{Here we compute the zero level-sets for the HJE arising from the quadcopter. The $x$-axis is the $x_1$-position of the quadcopter, and the $y$-axis is the angular position in the $\psi_1$ coordinate. The zero level-sets are computed for $t=0.025$, $0.05$, and $0.075$. This is a 12-dimensional, nonlinear optimal control problem.}
		\end{figure}

\begin{figure}[]
		\centering
	    \label{fig:quadcoptercharcurves}\includegraphics[scale=0.5]{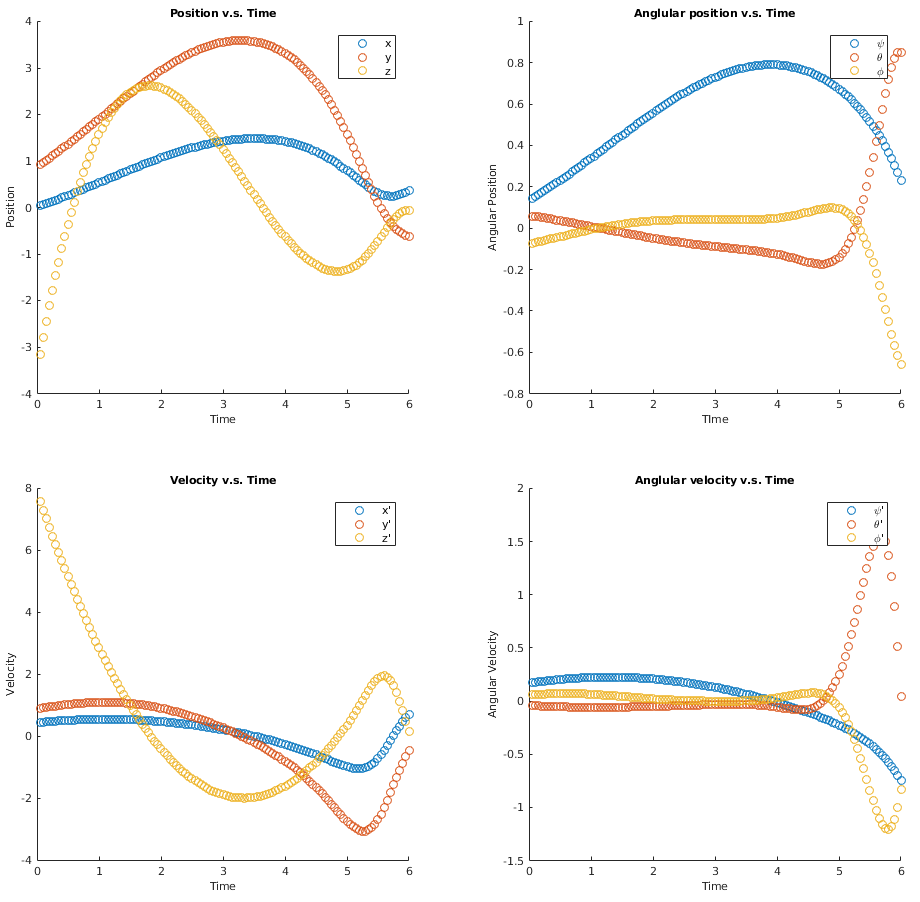}
	    \caption{Here we compute the characteristic curves/optimal trajectories for the quadcopter. We are computing at the terminal point $\tbf{x} = (0.36, -0.62, -0.06, 0.23, 0.85, -0.66, 0.72, -0.45, 0.15, -0.75, 0.04, -0.83)$ and we are computing at the terminal time $t = 6$ seconds. A plot of the trajectories computing using a different algorithm -- SQP -- looks identical.}
		\end{figure}

\begin{figure}[]
		\centering
	    \label{fig:quadcopter3D}\includegraphics[scale=0.5]{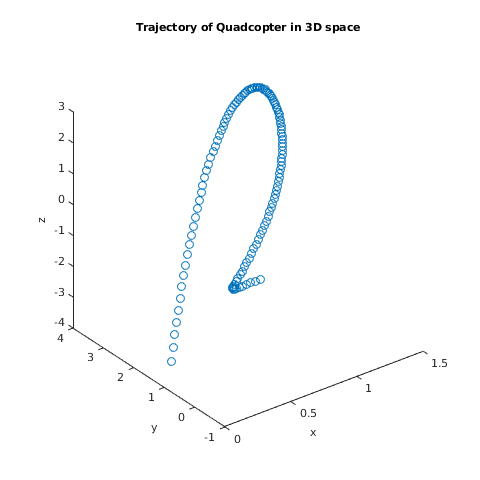}
	    \caption{We plot the $(x,y,z)$ coordinates of the quadcopter to give a plot of the trajectory of the quadcopter in 3D space.}
		\end{figure}

	\cref{fig:quadcopterhje} gives the zero level-sets of the HJE, and \cref{fig:quadcoptercharcurves} gives the result of computing the curves/trajectories. \cref{fig:quadcopter3D} plots the $x$, $y$, and $z$ positions of the quadcopter as it moves through time.
	
\newpage
\section{Discussion and Future Work}
\label{sec:disc}

	In this paper, we have presented a splitting method to compute solutions to general (i.e. convex and nonconvex) Hamilton-Jacobi equations which arise from general (i.e. linear and nonlinear) optimal control problems and general differential games problems.
	
	Some nice properties of our algorithm include: (1) relatively fast computation of solutions in high-dimensions, especially when we parallelize the algorithm which is embarrassingly parallelizable \cite{Herlihy:2008:AMP:1734069} (2) it can generate optimal trajectories of the optimal control/differential games problems, (3) it can compute problems with non-linear ODEs, (4) it can compute solutions for nonconvex Hamiltonians, (5) and the algorithm is embarrassingly parallelizable, i.e. each core can use the algorithm to compute the solution at a point, so given $N$ cores we can compute solutions of the HJE at $N$ points simultaneously.
	
	Splitting applied to optimal control problems has been used by \cite{6422363} where they apply it to cost functionals having a quadratic and convex term. In terms of Hamilton-Jacobi equations, the authors in Kirchner et al. \cite{KHDO18} (2018) effectively applied it to Hamilton-Jacobi equations arising from linear optimal control problems by using the Hopf formula. They make use of the Hopf formula and the closed-form solution to linear ODEs to not only solve HJE, but to also directly compute optimal trajectories in high-dimensional systems. The authors of this current paper have been working in parallel and also applied splitting to HJE and trajectory generation for nonlinear optimal control problems and minimax differential games.
	
	On a related note, see also previous work by Kirchner et al. \cite{8231132} where they apply the Hopf formula to differential games and show that complex ``teaming" behavior can arise, even with linearized pursuit-evasion models.
	
	As far as we know, the idea to use splitting for differential games problem for the discretization in equation \eqref{eq:discgenlaxhopfdg} and \eqref{eq:discgenlaxhopfdghopf} is new. And we believe it is worth examining if this PDHG-inspired method to solve minimax/saddle-point problems may apply to more general minimax/saddle-point optimization problems.

	The proof of convergence and the proof of approximation for our algorithm is still a work-in-progress. But it seems to be that for the examples in \cref{sec:examples}, we get relatively the correct answer, and our algorithm seems to scale linearly with dimension for even a nonlinear optimal control problem requiring nonconvex optimization (see \cref{subsubsec:dimscale}).
	
	We also believe that due to the deep connection between Hamilton-Jacobi equations and optimization methods (\cref{subsec:remarkhjeopt}), it is worthwhile to examine why our algorithm works. Not only that, but for our examples in \cref{sec:examples}, we were able to perform nonconvex optimization with only a single initial guess, whereas coordinate descent required multiple. And the authors also believe it is worth examining the algorithms \cref{alg:algohjedglax} and \cref{alg:algohjedghopf} as the computation of minimax differential games problems using a splitting method seems new. In essence, it may be possible to generalize these algorithms to apply to general minimax/saddle-point problems with continuous constraints.
	
	Some improvements to our algorithm for differential games problems (\cref{alg:algohjedglax} and \cref{alg:algohjedghopf}) can be foreseen: 
		\begin{enumerate}
		\item We have found speed-ups to our algorithm when we use acceleration methods \cite{Chambolle2011, Chambolle2016}
		\item One may be able to devise a more sophisticated stopping criteria as that in Kirchner et al. \cite{KHDO18}, where they apply a step-size-dependent stopping criteria based off work by \cite{goldstein2015adaptive}.
		\item We would also like to utilize higher-order approximations for the ODE and integral when discretizing the value functions of the optimal control or differential games problems. We note that for the Lax discretizations (\cref{alg:algohjeoclax} and \cref{alg:algohjedglax}), one can average the forward and backward Euler approximations to obtain higher accuracy, analogous to how the trapezoidal approximation is the average of the two.
		\item And we believe we might be able to make use of having a closed-form solution, or an approximate solution, to computing the characteristic curves, i.e. closed form solutions to $\frac{d}{dt}\tbf{x}(s) = H_p(\tbf{x}(s), \tbf{p}(s), s)$ and $\frac{d}{dt}\tbf{p}(s) = -H_x(\tbf{x}(s), \tbf{p}(s), s)$, much as in \cite{KHDO18}, where they make use of having a closed-form solution to linear differential equations by utilizing the exponential operator.
		\item There could be an advantage in combining the splitting method to pseudo-spectral methods \cite{10.1007/978-3-540-45056-6_21}.
		\item In the algorithms for differential equations \cref{alg:algohjedglax} and \cref{alg:algohjedghopf}, we are solving a saddle-point problem using gradients. We might obtain faster convergence if we used a Hessian-inspired method, such as split form of BFGS.
		\end{enumerate}			

\section{Acknowledgements}

	We give our deepest thanks to Matthew R. Kirchner and Gary Hewer for their enlightening discussions and suggested edits during this work. They provided us with a wealth of information on the history of the subject and existing methods, as well as important problems. And they were generous in giving suggestions on future directions.
	

\appendix

\newpage
\section{A practical tutorial for implementation purposes}
\label{appa:sec:practut}

\subsection{Optimal Control}

	Suppose we want to compute the Hamilton-Jacobi equations associated to the following optimal control problem:
		\begin{align}\label{eq:app:valphi}
		\phi(x,t) = \min_{\textbf{x}(\cdot), \tbf{u}(\cdot)} \l\{ g(\tbf{x}(0)) + \int_0^t L(\textbf{x}(s), \tbf{u}(s), s)\,ds \r\}
		\end{align}
	where $\tbf{x}(\cdot)$ and $\tbf{u}(s)$ satisfy the ODE
		\begin{align}\label{eq:app:ode}
		\l\{\begin{array}{l}
		\dot{\tbf{x}}(s) = \tbf{f}(\tbf{x}(s), \tbf{u}(s), s), \quad 0 < s < t \\
		\tbf{x}(t) = x
		\end{array}\r.
		\end{align}
	Here $(x,t)\in\mbb{R}^n\times [0, \infty)$ are fixed, and is the point that we want to compute the HJE solution $\varphi$.
	
	Then we can use \cref{alg:algohjeoclax} or \cref{alg:algohjeochopf}, which we describe in more detail below.

\subsubsection{Practical tutorial for \cref{alg:algohjeoclax}}
	If we want to use the discretized Lax formula (with a backward Euler discretization of the ODE dynamics), then:

\begin{enumerate}	

	\item Discretize the time domain:
		\begin{align*}
		0 = s_0 < s_1 < s_2 < \cdots < s_{N-1} < s_N = t.
		\end{align*}
	In our numerical experiments, we chose a uniform discretization of size $\delta \defeq t/N$.
	
	\item Approximate \eqref{eq:app:ode} using backward Euler, and also discretize \eqref{eq:app:valphi} to obtain (as in \eqref{eq:disclaxBE}),
		\begin{align*}
		\phi(x,t) \approx \max_{\{p_j\}_{j=1}^{N}} \min_{\{x_j\}_{j=0}^{N}} \l\{ g(x_0) + \sum_{j=1}^{N} \l<p_{j}, x_{j} - x_{j-1}\r>  - \delta \sum_{j=1}^{N} H(x_j, p_j, s_j) \r\}
		\end{align*}
	where $x_j = \tbf{x}(s_j)$, and $p_j = \tbf{p}(x_j)$. Let us denote $x = x_{\text{target}}$ to clarify notation. 

	\item Initialize: 
		\begin{enumerate}
			\item Choose $\delta > 0$, set $N = t/\delta$ (although note that since we are using the zero-th time-step, then we are updating $N+1$ points).
	
			\item Randomly initialize $\tilde{x}^0 \defeq (x^0_0, x^0_1, \ldots, x^0_{N-1}, x^0_N)$, but with $x_N^0 \equiv x_{\text{target}}$.
	
			\item Randomly initialize $\tilde{p}^0 \defeq (p^0_0, p^0_1, \ldots, p^0_{N-1}, p^0_N)$, but with $p_0^0 \equiv 0$, as we won't be updating $p_0^0$; it is only there for computational accounting.
	
			\item Set $\tilde{z}^0 \defeq (z_0^0, z_1^0, \ldots, z_N^0) = (x^0_0, x^0_1, \ldots, x^0_{N-1}, x^0_N)$.
	
			\item Choose $\sigma, \tau$ such that $\sigma \tau < 1/\|D\|_2^2 = 0.25$ and $\theta \in [0,1]$ (we suggest $\theta = 1$).
			
			\item Choose some tolerance $\text{tol} > 0$ small.
		\end{enumerate}
	
	\item Set
		\begin{align*}
		D = \begin{pmatrix}
		\tbf{0} & \tbf{0} & \tbf{0} & \tbf{0} & \cdots & \tbf{0} & \tbf{0} \\
		-I & I & \tbf{0} & \tbf{0} & \cdots & \tbf{0} & \tbf{0} \\
		\tbf{0} & -I & I & \tbf{0} & \cdots & \tbf{0} & \tbf{0} \\
		\vdots & \vdots & \vdots & \vdots & \ddots & \vdots & \vdots\\
		\tbf{0} & \tbf{0} & \tbf{0} & \tbf{0} & \cdots & \tbf{0} & I \\
		\tbf{0} & \tbf{0} & \tbf{0} & \tbf{0} & \cdots & \tbf{0} & -I
		\end{pmatrix}
		\end{align*}
	where $\tbf{0}$ is a $(\text{dim}\times\text{dim})$ zero matrix, where $\text{dim}$ is the size of the space variable, and $I$ is the $(\text{dim}\times\text{dim})$ identity matrix. 
	
	Note that in the algorithm, we can replace $(D\tilde{z}^k)_j = z^k_j - z^k_{j-1}$, and similarly with $(D^T \tilde{p}^k)_j = p^k_j - p^k_{j+1}$, so we can save time by not performing a full matrix multiplication.
	
	\item Then perform the algorithm found in \cref{alg:appa:algooclax}.

\end{enumerate}

	\begin{algorithm}
		\caption{Practical tutorial for the Lax formula with backward Euler, for Optimal Control}
		\begin{algorithmic}\label{alg:appa:algooclax}
		\STATE{Given: $x_{\text{target}}\in\mbb{R}^d$ and time $t\in (0, \infty)$.}
		\STATE{\vspace{1mm}}
		\WHILE{$(\|\tilde{x}^{k+1} - \tilde{x}^k\|_2^2 > \text{tol} \text{ or } \|\tilde{p}^{k+1} - \tilde{p}^k\|^2_2 > \text{tol}) \text{ and } (\text{count} < \text{max\_count})$ }
		\FOR{$j = 1 \txt{ to } N$}
		\STATE{${p_j^{k+1}} = \argmin_{p} \l\{ \delta H(x_j^k, p, s_j) + \frac{1}{2\sigma}\|p - (p_j^k + \sigma(D\tilde{z}^k)_j) \|_2^2\r\}$}
		\ENDFOR
		\STATE{\vspace{1mm}}
		\FOR{$j = 0$}
		\STATE{$x_0^{k+1} = \argmin_{x} \l\{ g(x) + \frac{1}{2\tau}\|x - (x_0^{k} - \tau(D^T\tilde{p}^k)_0)\|_2^2 \r\} \quad \text{ (note $p_0^k = 0$)}
			$}
		\ENDFOR
		\FOR{$j = 1 \txt{ to } N-1$}
		\STATE{$x_j^{k+1} = \argmin_{x} \l\{ - \delta H(x, p_j^{k+1}, s_j) + \frac{1}{2\tau}\|x - (x_j^{k} - \tau(D^T\tilde{p}^k)_j)\|_2^2
			\r\}$}
		\ENDFOR
		\STATE{\vspace{1mm}}
		\FOR{$j = 0 \txt{ to } N$}
		\STATE{$z_j^{k+1} = x_j^{k+1} + \theta (x_j^{k+1} - x_j^k)$}
		\ENDFOR
		\ENDWHILE
		\STATE{\vspace{1mm}}
		\STATE{$\text{fval} = g(x_0) + \sum_{j=1}^{N} \l<p_{j}, x_{j} - x_{j-1}\r>  - \delta \sum_{j=1}^{N} H(x_j, p_j, s_j)$}
		\RETURN fval
		\end{algorithmic}
	\end{algorithm}

\newpage
\subsubsection{Practical tutorial for \cref{alg:algohjeochopf}}

	If we want to use the discretized Hopf formula, then:
	\begin{enumerate}
		\item Discretize the time domain:
		\begin{align*}
		0 = s_0 < s_1 < s_2 < \cdots < s_{N-1} < s_N = t.
		\end{align*}
	In our numerical experiments, we chose a uniform discretization of size $\delta \defeq t/N$.
	
	\item Approximate \eqref{eq:app:ode} using backward Euler, and also discretize \eqref{eq:app:valphi} to obtain (as in \eqref{eq:dischopf}),
		\begin{align*}
		\phi(x,t) \approx \max_{\{p_j\}_{j=1}^N} \min_{\{x_j\}_{j=1}^N } \l\{ -g^*(p_1) + \l<p_N, x\r> + \sum_{j=1}^{N-1} \l<p_j - p_{j+1}, x_j\r>  - \delta \sum_{j=1}^{N} H(x_{j}, p_{j}, s_{j})\r\} 
		\end{align*}
	where $x_j = \tbf{x}(s_j)$, and $p_j = \tbf{p}(x_j)$. Let us denote $x = x_{\text{target}}$ to clarify notation. 

	\item Initialize: 
		\begin{enumerate}
			\item Choose $\delta > 0$, set $N = t/\delta$.
	
			\item Randomly initialize 	$\tilde{x}^0 \defeq (x^0_1, \ldots, x^0_{N-1}, x^0_N)$, but with $x_N^0 \equiv x_{\text{target}}$.
	
			\item Randomly initialize $\tilde{p}^0 \defeq (p^0_1, \ldots, p^0_{N-1}, p^0_N)$.
	
			\item Set $\tilde{z}^0 \defeq (z_1^0, \ldots, z_N^0) = (x^0_1, \ldots, x^0_{N-1}, x^0_N)$.
	
			\item Choose $\sigma, \tau$ such that $\sigma \tau < 1/\|D\|_2^2 = 0.25$ and $\theta \in [0,1]$ (we suggest $\theta = 1$).
			
			\item Choose some tolerance $\text{tol}>0$ small.
		\end{enumerate}
		
	\item Set
		\begin{align*}
		D = \begin{pmatrix}
		I & -I & \tbf{0} & \tbf{0} & \cdots & \tbf{0} & \tbf{0} \\
		\tbf{0} & I & -I & \tbf{0} & \cdots & \tbf{0} & \tbf{0} \\
		\vdots & \vdots & \vdots & \vdots & \ddots & \vdots & \vdots\\
		\tbf{0} & \tbf{0} & \tbf{0} & \tbf{0} & \cdots & \tbf{0} & -I \\
		\tbf{0} & \tbf{0} & \tbf{0} & \tbf{0} & \cdots & \tbf{0} & I
		\end{pmatrix}
		\end{align*}
	where $\tbf{0}$ is a $(\text{dim}\times\text{dim})$ zero matrix, where $\text{dim}$ is the size of the space variable, and $I$ is the $(\text{dim}\times\text{dim})$ identity matrix. 
	
	Note we can replace $(D^T\tilde{z^k})_j = z^k_j - z^k_{j-1}$ and $(D \tilde{p}^k)_j = p^k_{j} - p^k_{j+1}$, so we can save time by not performing a full matrix multiplication. 
	
	Also note that the $D$ here is different than in \cref{alg:algohjeoclax} and \cref{alg:appa:algooclax}.
	
	\item Then perform the algorithm found in \cref{alg:appa:algoochopf}
	
	\end{enumerate}

	\begin{algorithm}
		\caption{Practical tutorial for the Hopf formula, for Optimal Control}
		\begin{algorithmic}\label{alg:appa:algoochopf}
		\STATE{Given: $x_{\text{target}}\in\mbb{R}^d$ and time $t\in (0, \infty)$.}
		\STATE{\vspace{1mm}}
		\WHILE{$(\|\tilde{x}^{k+1} - \tilde{x}^k\|_2^2 > \text{tol} \text{ or } \|\tilde{p}^{k+1} - \tilde{p}^k\|^2_2 > \text{tol}) \text{ and } (\text{count} < \text{max\_count})$ }
		\FOR{$j = 1$}
		\STATE{${p_j^{k+1}} = \argmin_{p} \l\{ g^*(p) + \delta H(x_1^k, p, s_1) + \frac{1}{2\sigma}\|p - (p_1^k + \sigma(D^T\tilde{z}^k)_1) \|_2^2\r\}$}
		\ENDFOR
		\FOR{$j = 2 \txt{ to } N$}
		\STATE{${p_j^{k+1}} = \argmin_{p} \l\{ \delta H(x_j^k, p, s_j) + \frac{1}{2\sigma}\|p - (p_j^k + \sigma(D^T\tilde{z}^k)_j) \|_2^2\r\}$}
		\ENDFOR
		\STATE{\vspace{1mm}}
		\FOR{$j = 1 \txt{ to } N-1$}
		\STATE{$x_j^{k+1} = \argmin_{x} \l\{ - \delta H(x, p_j^{k+1}, s_j) + \frac{1}{2\tau}\|x - (x_j^{k} - \tau(D\tilde{p}^k)_j)\|_2^2
			\r\}$}
		\ENDFOR
		\STATE{\vspace{1mm}}
		\FOR{$j = 1 \txt{ to } N$}
		\STATE{$z_j^{k+1} = x_j^{k+1} + \theta (x_j^{k+1} - x_j^k)$}
		\ENDFOR
		\ENDWHILE
		\STATE{\vspace{1mm}}
		\STATE{$\text{fval} = -g^*(p_1) + \l<p_N, x_{\text{target}}\r> + \sum_{j=1}^{N-1} \l<p_j - p_{j+1}, x_j\r>  - \delta \sum_{j=1}^{N} H(x_{j}, p_{j}, s_{j})$}
		\RETURN fval
		\end{algorithmic}
	\end{algorithm}

\newpage
\subsection{Differential Games}

	Suppose we want to solve the differential games problem with the following dynamics,
		\begin{align}\label{eq:appa:dgode}
		\l\{ \begin{array}{l}
		\begin{pmatrix}\dot{\tbf{x}}(s) \\ \dot{\tbf{y}}(s)) \end{pmatrix}= \begin{pmatrix}\tbf{f}_1(\tbf{x}(s), \tbf{y}(s), \bsym{\alpha}(s), \bsym{\beta}(s), s) \\ \tbf{f}_2(\tbf{x}(s), \tbf{y}(s), \bsym{\alpha}(s), \bsym{\beta}(s), s) \end{pmatrix} \quad 0<s<t \\
		\begin{pmatrix} \tbf{x}(t) \\ \tbf{y}(t) \end{pmatrix} = \begin{pmatrix} x \\ y \end{pmatrix}
		\end{array}\r.
		\end{align}
	and with the following value function,
		\begin{align}\label{eq:appa:valdg}
		\phi(x,y,t) = \inf_{\bsym{\alpha}(\cdot), \tbf{x}(\cdot)} \sup_{\bsym{\beta}(\cdot), \tbf{y}(\cdot)} \l\{ g(\tbf{x}(0), \tbf{y}(0)) + \int_0^t L(\tbf{x}(s), \tbf{y}(s), \bsym{\alpha}(s), \bsym{\beta}(s), s)\,ds\r\}
		\end{align}
	Then we can discretize the above equation and use \cref{alg:algohjedglax} and \cref{alg:algohjedghopf}, which we describe in more detail below.
	
\subsubsection{Practical tutorial for \cref{alg:algohjedglax}}

	If we want to use the discretized Lax formula for differential games (with a backward Euler discretization of the ODE dynamics) \eqref{eq:discgenlaxhopfdg}, then:
	
	\begin{enumerate}	

	\item Discretize the time domain:
		\begin{align*}
		0 = s_0 < s_1 < s_2 < \cdots < s_{N-1} < s_N = t.
		\end{align*}
	In our numerical experiments, we chose a uniform discretization of size $\delta \defeq t/N$.
	
	\item Approximate \eqref{eq:appa:dgode} using backward Euler, and also discretize \eqref{eq:appa:valdg} to obtain (as in \eqref{eq:discgenlaxhopfdg}),
		\begin{align*}
		\hspace{-1in}\phi(x,y,t) \approx \min_{\{q_j\}_{j=1}^N} \max_{\{p_j\}_{j=1}^N} \min_{\{x_j\}_{j=0}^N} \max_{\{y_j\}_{j=0}^N} \l\{ g(x_0, y_0) + \sum_{j=1}^N \l<\begin{pmatrix}p_j \\ -q_j \end{pmatrix}, \begin{pmatrix}x_j - x_{j-1} \\ y_j - y_{j-1} \end{pmatrix}\r> - \delta \sum_{j=1}^N H(x_j, y_j, p_j, -q_j, s_j)  \r\}
		\end{align*}
	where $x_j = \tbf{x}(s_j)$, and similarly for $y_j$, $q_j$, and $p_j$. Let us denote $x = x_{\text{target}}$ and $y = y_{\text{target}}$ to clarify notation. 

	\item Initialize: 
		\begin{enumerate}[1.]
			\item Choose $\delta > 0$, set $N = t/\delta$ (although note that since we are using the zero-th time-step, then we are updating $N+1$ points).
	
			\item Randomly initialize $\tilde{x}^0 \defeq (x^0_0, x^0_1, \ldots, x^0_{N-1}, x^0_N)$, but with $x_N^0 \equiv x_{\text{target}}$, and similarly for $\tilde{y}^0$.
	
			\item Randomly initialize $\tilde{p}^0 \defeq (p^0_0, p^0_1, \ldots, p^0_{N-1}, p^0_N)$, but with $p_0^0 \equiv 0$, as we won't be updating $p_0^0$; it is only there for computational accounting. Do a similary initialization for $\tilde{q}^0$.
	
			\item Set $\tilde{z}^0 \defeq (z_0^0, z_1^0, \ldots, z_N^0) = (x^0_0, x^0_1, \ldots, x^0_{N-1}, x^0_N)$, and set $\tilde{w}^0 \defeq (w_0^0, w_1^0, \ldots, w_N^0) = (y^0_0, y^0_1, \ldots, y^0_{N-1}, y^0_N)$.
	
			\item Choose $\sigma, \tau$ such that $\sigma \tau < 1/\|D\|_2^2 = 0.25$ and $\theta \in [0,1]$ (we suggest $\theta = 1$).
			
			\item Choose some tolerance $\text{tol} > 0$ small.
		\end{enumerate}
	
	\item Set
		\begin{align*}
		D = \begin{pmatrix}
		\tbf{0} & \tbf{0} & \tbf{0} & \tbf{0} & \cdots & \tbf{0} & \tbf{0} \\
		-I & I & \tbf{0} & \tbf{0} & \cdots & \tbf{0} & \tbf{0} \\
		\tbf{0} & -I & I & \tbf{0} & \cdots & \tbf{0} & \tbf{0} \\
		\vdots & \vdots & \vdots & \vdots & \ddots & \vdots & \vdots\\
		\tbf{0} & \tbf{0} & \tbf{0} & \tbf{0} & \cdots & \tbf{0} & I \\
		\tbf{0} & \tbf{0} & \tbf{0} & \tbf{0} & \cdots & \tbf{0} & -I
		\end{pmatrix}
		\end{align*}
	where $\tbf{0}$ is a $(\text{dim}\times\text{dim})$ zero matrix, where $\text{dim}$ is the size of the space variable, and $I$ is the $(\text{dim}\times\text{dim})$ identity matrix. Note this is a very sparse matrix and we take advantage of this.
	
	Note that we can replace $D\tilde{z}^k = z^k_j - z^k_{j-1}$, and $D\tilde{w}^k = w^k_j - w^k_{j-1}$, and $D^T\tilde{p}^k = p^k_j - p^k_{j+1}$, and $D^T\tilde{q}^k = q^k_j - q^k_{j+1}$, so we don't have to perform the full matrix multiplication.
	
	\item Then perform the algorithm found in \cref{alg:appa:algodglax}.
	
\end{enumerate}
	
	\begin{algorithm}
		\caption{Practical tutorial for the Lax formula with backward Euler, for Differential Games}
		\begin{algorithmic}\label{alg:appa:algodglax}
		\STATE{Given: $x_{\text{target}}\in\mbb{R}^d$ and time $t\in (0, \infty)$.}
		\STATE{\vspace{1mm}}
		\WHILE{$(\|\tilde{x}^{k+1} - \tilde{x}^k\|_2^2 > \text{tol} \text{ or } \|\tilde{p}^{k+1} - \tilde{p}^k\|^2_2 > \text{tol} \text{ or } \|\tilde{y}^{k+1} - \tilde{y}^k\|_2^2 > \text{tol} \text{ or } \|\tilde{q}^{k+1} - \tilde{q}^k\|^2_2 > \text{tol}) \text{ and } (\text{count} < \text{max\_count})$ }
		\FOR{$j = 1 \txt{ to } N$}
		\STATE{${p_j^{k+1}} = \argmin_{p} \l\{ \delta H(x_j^k, y_j^k, p, -q_j^k, s_j) + \frac{1}{2\sigma}\|p - (p_j^k + \sigma(D\tilde{z}^k)_j) \|_2^2\r\}$}
		\ENDFOR
		\STATE{\vspace{1mm}}
		
		\FOR{$j=1 \text{ to } N$}
		\STATE{${q_j^{k+1}} = \argmin_{q} \l\{ -\delta H(x_j^k, y_j^k, p_j^{k+1}, -q, s_j) + \frac{1}{2\sigma}\|q - (q_j^k + \sigma(D\tilde{w}^k)_j) \|_2^2\r\}$}
		\ENDFOR
		\STATE{\vspace{1mm}}

		\FOR{$j = 0$}
		\STATE{$x_0^{k+1} = \argmin_{x} \l\{ g(x, y_0^{k}) + \frac{1}{2\tau}\|x - (x_0^{k} - \tau (D^T\tilde{p}^k)_0)\|_2^2
			\r\}$}
		\ENDFOR
		\FOR{$j = 1 \txt{ to } N-1$}
		\STATE{$x_j^{k+1} = \argmin_{x} \l\{ - \delta H(x, y_j^k, p_j^{k+1}, -q_j^{k+1}, s_j) + \frac{1}{2\tau}\|x - (x_j^{k} - \tau (D^T\tilde{p}^k)_j)\|_2^2
			\r\}$}
		\ENDFOR
		\STATE{\vspace{1mm}}
		
		\FOR{$j = 0$}
		\STATE{$y_0^{k+1} = \argmin_{y} \l\{ -g(x_0^{k+1}, y) + \frac{1}{2\tau}\|y - (y_0^{k} - \tau (D^T\tilde{q}^k)_0)\|_2^2
			\r\}$}
		\ENDFOR
		\FOR{$j = 1 \txt{ to } N-1$}
		\STATE{$y_j^{k+1} = \argmin_{y} \l\{ - \delta H(x_j^{k+1}, y, p_j^{k+1}, -q_j^{k+1}, s_j) + \frac{1}{2\tau}\|y - (y_j^{k} - \tau (D^T\tilde{q}^k)_j)\|_2^2
			\r\}$}
		\ENDFOR		
		\STATE{\vspace{1mm}}
		
		\FOR{$j = 0 \txt{ to } N$}
		\STATE{$z_j^{k+1} = x_j^{k+1} + \theta (x_j^{k+1} - x_j^k)$}
		\STATE{$w_j^{k+1} = y_j^{k+1} + \theta (y_j^{k+1} - y_j^k)$}
		\ENDFOR
		
		\ENDWHILE
		\STATE{\vspace{1mm}}
		
		\STATE{$\text{fval} = g(x_0, y_0) + \sum_{j=1}^N \l<\begin{pmatrix}p_j \\ -q_j \end{pmatrix}, \begin{pmatrix}x_j - x_{j-1} \\ y_j - y_{j-1} \end{pmatrix}\r> - \delta \sum_{j=1}^N H(x_j, y_j, p_j, -q_j, s_j)$}
		\RETURN fval
		\end{algorithmic}
	\end{algorithm}
	
\newpage
\subsubsection{Practical tutorial for \cref{alg:algohjedghopf}}

	If we want to use the discretized Hopf formula for differential games (with a backward Euler discretization of the ODE dynamics) \eqref{eq:discgenlaxhopfdghopf}, then:
	
	\begin{enumerate}	

	\item Discretize the time domain:
		\begin{align*}
		0 = s_0 < s_1 < s_2 < \cdots < s_{N-1} < s_N = t.
		\end{align*}
	In our numerical experiments, we chose a uniform discretization of size $\delta \defeq t/N$.
	
	\item Approximate \eqref{eq:appa:dgode} using backward Euler, and also discretize \eqref{eq:appa:valdg} to obtain (as in \eqref{eq:discgenlaxhopfdghopf}),
		\begin{equation*}
		\begin{aligned}
		\hspace{-1.1in}\phi(x,y,t) &\approx \min_{\{q_j\}_{j=1}^N} \max_{\{p_j\}_{j=1}^N} \min_{\{x_j\}_{j=1}^N} \max_{\{y_j\}_{j=1}^N} \l\{ -e^*(p_1) - h_*(-q_1) + \l<\begin{pmatrix} p_N \\ -q_N \end{pmatrix}, \begin{pmatrix} x \\ y \end{pmatrix} \r> + \sum_{j=1}^{N-1} \l<\begin{pmatrix}p_j - p_{j+1} \\ -(q_j - q_{j+1}) \end{pmatrix}, \begin{pmatrix}x_j \\ y_j \end{pmatrix}\r>\r. \\
		& \hspace{4.5in}\l. - \delta \sum_{j=1}^N H(x_j, y_j, p_j, -q_j, s_j)  \r\}
		\end{aligned}
		\end{equation*}
	where $x_j = \tbf{x}(s_j)$, and similarly for $y_j$, $q_j$, and $p_j$. Let us denote $x = x_{\text{target}}$ and $y = y_{\text{target}}$ to clarify notation. 

	\item Initialize: 
		\begin{enumerate}[1.]
			\item Choose $\delta > 0$, set $N = t/\delta$.
	
			\item Randomly initialize $\tilde{x}^0 \defeq (x^0_0, x^0_1, \ldots, x^0_{N-1}, x^0_N)$, but with $x_N^0 \equiv x_{\text{target}}$, and similarly for $\tilde{y}^0$.
	
			\item Randomly initialize $\tilde{p}^0 \defeq (p^0_0, p^0_1, \ldots, p^0_{N-1}, p^0_N)$, but with $p_0^0 \equiv 0$, as we won't be updating $p_0^0$; it is only there for computational accounting. Do the same initialization for $\tilde{q}^0$.
	
			\item Set $\tilde{z}^0 \defeq (z_0^0, z_1^0, \ldots, z_N^0) = (x^0_0, x^0_1, \ldots, x^0_{N-1}, x^0_N)$, and set $\tilde{w}^0 \defeq (w_0^0, w_1^0, \ldots, w_N^0) = (y^0_0, y^0_1, \ldots, y^0_{N-1}, y^0_N)$.
	
			\item Choose $\sigma, \tau$ such that $\sigma \tau < 1/\|D\|_2^2 = 0.25$ and $\theta \in [0,1]$ (we suggest $\theta = 1$).
			
			\item Choose some tolerance $\text{tol} > 0$ small.
		\end{enumerate}
	
	\item Set
		\begin{align*}
		D = \begin{pmatrix}
		I & -I & \tbf{0} & \tbf{0} & \cdots & \tbf{0} & \tbf{0} \\
		\tbf{0} & I & -I & \tbf{0} & \cdots & \tbf{0} & \tbf{0} \\
		\vdots & \vdots & \vdots & \vdots & \ddots & \vdots & \vdots\\
		\tbf{0} & \tbf{0} & \tbf{0} & \tbf{0} & \cdots & \tbf{0} & -I \\
		\tbf{0} & \tbf{0} & \tbf{0} & \tbf{0} & \cdots & \tbf{0} & I
		\end{pmatrix}
		\end{align*}
	where $\tbf{0}$ is a $(\text{dim}\times\text{dim})$ zero matrix, where $\text{dim}$ is the size of the space variable, and $I$ is the $(\text{dim}\times\text{dim})$ identity matrix. Note this is a very sparse matrix and we take advantage of this.
	
	Alose note that we can replace $D^T\tilde{z}^k = z^k_j - z^k_{j-1}$, and $D^T\tilde{w}^k = w^k_j - w^k_{j-1}$, and $D\tilde{p}^k = p^k_j - p^k_{j+1}$, and $D\tilde{q}^k = q^k_j - q^k_{j+1}$, so we don't have to perform the full matrix multiplication.
	
	\item Then perform the algorithm found in \cref{alg:appa:algodghopf}.
	
\end{enumerate}
	
	\begin{algorithm}
		\caption{Practical tutorial for the Lax formula with backward Euler, for Differential Games}
		\begin{algorithmic}\label{alg:appa:algodghopf}
		\STATE{Given: $x_{\text{target}}\in\mbb{R}^d$ and time $t\in (0, \infty)$.}
		\STATE{\vspace{1mm}}
		\WHILE{$(\|\tilde{x}^{k+1} - \tilde{x}^k\|_2^2 > \text{tol} \text{ or } \|\tilde{p}^{k+1} - \tilde{p}^k\|^2_2 > \text{tol} \text{ or } \|\tilde{y}^{k+1} - \tilde{y}^k\|_2^2 > \text{tol} \text{ or } \|\tilde{q}^{k+1} - \tilde{q}^k\|^2_2 > \text{tol}) \text{ and } (\text{count} < \text{max\_count})$ }
		\FOR{$j = 1$}
		\STATE{${p_1^{k+1}} = \argmin_{p} \l\{ e^*(p) + \delta H(x_1^k, y_1^k, p, -q_1^k, s_1) + \frac{1}{2\sigma}\|p - (p_1^k + \sigma(D^T\tilde{z}^k)_1) \|_2^2\r\}$}
		\ENDFOR
		\FOR{$j = 2 \txt{ to } N$}
		\STATE{${p_j^{k+1}} = \argmin_{p} \l\{ \delta H(x_j^k, y_j^k, p, -q_j^k, s_j) + \frac{1}{2\sigma}\|p - (p_j^k + \sigma(D^T\tilde{z}^k)_j) \|_2^2\r\}$}
		\ENDFOR
		\STATE{\vspace{1mm}}
		
		\FOR{$j = 1$}
		\STATE{${q_1^{k+1}} = \argmin_{q} \l\{- h_*(-q) -\delta H(x_1^k, y_1^k, p_1^{k+1}, -q, s_1) + \frac{1}{2\sigma}\|q - (q_1^k + \sigma(D^T\tilde{w}^k)_1) \|_2^2\r\}$}
		\ENDFOR
		\FOR{$j=2 \text{ to } N$}
		\STATE{${q_j^{k+1}} = \argmin_{q} \l\{ -\delta H(x_j^k, y_j^k, p_j^{k+1}, -q, s_j) + \frac{1}{2\sigma}\|q - (q_j^k + \sigma(D^T\tilde{w}^k)_j) \|_2^2\r\}$}
		\ENDFOR
		\STATE{\vspace{1mm}}

		\FOR{$j = 1 \txt{ to } N-1$}
		\STATE{$x_j^{k+1} = \argmin_{x} \l\{ - \delta H(x, y_j^k, p_j^{k+1}, -q_j^{k+1}, s_j) + \frac{1}{2\tau}\|x - (x_j^{k} - \tau (D\tilde{p}^k)_j)\|_2^2
			\r\}$}
		\ENDFOR
		\STATE{\vspace{1mm}}
		
		\FOR{$j = 1 \txt{ to } N-1$}
		\STATE{$y_j^{k+1} = \argmin_{y} \l\{ - \delta H(x_j^{k+1}, y, p_j^{k+1}, -q_j^{k+1}, s_j) + \frac{1}{2\tau}\|y - (y_j^{k} - \tau (D\tilde{q}^k)_j)\|_2^2
			\r\}$}
		\ENDFOR		
		\STATE{\vspace{1mm}}
		
		\FOR{$j = 1 \txt{ to } N$}
		\STATE{$z_j^{k+1} = x_j^{k+1} + \theta (x_j^{k+1} - x_j^k)$}
		\STATE{$w_j^{k+1} = y_j^{k+1} + \theta (y_j^{k+1} - y_j^k)$}
		\ENDFOR
		
		\ENDWHILE
		\STATE{\vspace{1mm}}
		
		\STATE{$\text{fval} = -e^*(p_1) - h_*(-q_1) + \l<\begin{pmatrix} p_N \\ -q_N \end{pmatrix}, \begin{pmatrix} x_{\text{target}} \\ y_{\text{target}} \end{pmatrix} \r> + \sum_{j=1}^{N-1} \l<\begin{pmatrix}p_j - p_{j+1} \\ -(q_j - q_{j+1}) \end{pmatrix}, \begin{pmatrix}x_j \\ y_j \end{pmatrix}\r> - \delta \sum_{j=1}^N H(x_j, y_j, p_j, -q_j, s_j) $}
		\RETURN fval
		\end{algorithmic}
	\end{algorithm}

\bibliographystyle{siamplain}
\bibliography{references}
\end{document}